\documentclass[mathpazo]{cicp}
\usepackage{amsmath,amssymb}
\usepackage[utf8]{inputenc}
\usepackage{mathrsfs}
\usepackage{todonotes}
\usepackage{color}

\newcommand{\gb}{\mathbf{g}}

\newcommand{\R}{\mathbb{R}}
\newcommand{\Vb}{\mathbf{V}}
\newcommand{\Vg}{\mathbf{V}_{\text{G}}}
\newcommand{\Vmfs}{\mathbf{V}_{\text{MFS}}}
\newcommand{\mi}{\mathrm{i}}
\newcommand{\varphib}{\boldsymbol{\varphi}}

\renewcommand{\Re}{\operatorname{Re}}

\newcommand{\uinc}{u^{\text{inc}}}
\newcommand{\utot}{u^{\text{tot}}}
\newcommand{\supp}{\operatorname{supp}}
\newcommand{\tstep}{{\Delta t}}
\renewcommand{\O}{\mathcal{O}}
\newcommand{\dist}{\operatorname{dist}}
\usepackage{graphicx}
\usepackage{subcaption}
\usepackage{graphicx}
\numberwithin{equation}{section}
\newcommand{\blue} [1] {#1}
\newcommand{\red} [1] {#1}

\begin{document}
%%%%% title : short title may not be used but TITLE is required.
% \title{TITLE}
% \title[short title]{TITLE}
\title{A modified convolution quadrature combined with the method of fundamental solutions and Galerkin BEM for acoustic scattering}
%==============
%Ebraheem Aldahham and Lehel Banjai\corrauth}
%\address{Mathematical and Computer Sciences, Heriot-Watt University, Edinburgh EH14 4AS, UK
%===============
%%%%% author(s) :
% single author:
% \author[name in running head]{AUTHOR\corrauth}
% [name in running head] is NOT OPTIONAL, it is a MUST.
% Use \corrauth to indicate the corresponding author.
% Use \email to provide email address of author.
% \footnote and \thanks are not used in the heading section.
% Another acknowlegments/support of grants, state in Acknowledgments section
% \section*{Acknowledgments}
%\author[O.~Author]{Only Author\corrauth}
%\address{School of Mathematical Sciences, Beijing Normal University,
%Beijing 100875, P.R. China}
%\email{{\tt author@email} (O.~Author)}

% multiple authors:
% Note the use of \affil and \affilnum to link names and addresses.
% The author for correspondence is marked by \corrauth.
% use \emails to provide email addresses of authors
% e.g. below example has 3 authors, first author is also the corresponding
%      author, author 1 and 3 having the same address.
%\author[Zhang Z R et.~al.]{Ebraheem Aldahham\affil{1}\comma\corrauth,
%       and Lehel Banjai\affil{1}}
\author[Aldahham E and Banjai L.]{Ebraheem Aldahham\affil{1}\comma\corrauth,
       and Lehel Banjai\affil{1}}
\address{\affilnum{1}\ Maxwell Institute for Mathematical Sciences, School of Mathematical and Computer Sciences, Heriot-Watt University, Edinburgh EH14 4AS, UK}
\emails{{\tt ea56@hw.ac.uk} (E.~Aldahham), {\tt l.banjai@hw.ac.uk} (L.~Banjai)}
 %\footnote and \thanks are not used in the heading section.
 %Another acknowlegments/support of grants, state in Acknowledgments section
 %\section*{Acknowledgments}

%%%%% Begin Abstract %%%%%%%%%%%
\begin{abstract}
We describe a numerical method for the solution of acoustic exterior scattering problems based on the time-domain boundary integral representation of the solution. As the spatial discretization of the resulting time-domain boundary integral equation we use either the method of fundamental solutions (MFS) or the Galerkin boundary element method (BEM). In time we apply either a standard convolution quadrature (CQ) based on an A-stable linear multistep method or a modified CQ scheme. It is well-known that the standard low-order CQ schemes for hyperbolic problems suffer from strong dissipation and dispersion properties. The modified scheme is designed to avoid these properties. We give a careful description of the modified scheme and its implementation with differences due to different spatial discretizations highlighted. Numerous numerical experiments illustrate the effectiveness of the modified scheme \blue{and dramatic improvement  with errors up to two orders of magnitude smaller in comparison with the standard scheme}.
\end{abstract}
%%%%% end %%%%%%%%%%%

%%%%% AMS/PACs/Keywords %%%%%%%%%%%
%\pac{}
%=====================
%I have found on the Internet the following
%https://mathscinet.ams.org/msnhtml/msc2020.pdf
%which related to ams
%=====================

%%%% maketitle %%%%%
\maketitle

%%%% Start %%%%%%
\section{Introduction}
In this paper we investigate a class of numerical methods for the scattering problem: Find $u(t) \in H^1(\Omega^+)$ such that
\begin{equation}\label{pde}
    \begin{aligned}
    \partial_t^2 u -\Delta u &= 0 & \text{for } (t,x) \in [0,T] \times \Omega^+\\
    u(0) = \partial_t u(0) &= 0 & \text{for } x \in \Omega^+\\
    u &= g & \text{for } (t,x) \in [0,T] \times \Gamma,
    \end{aligned}
\end{equation}
where $\Omega \subset \R^d$, $d = 2,3$,  is an open, bounded Lipschitz domain with exterior $\Omega^+ = \R^d \setminus \overline{\Omega}$ and boundary $\Gamma = \partial \Omega$. If $g$ is the trace of an incident wave $-\uinc$ on $\Gamma$ with $\supp \uinc(0) \subset \Omega^+$ then $u$ is the {\em scattered wave} and $\utot = u+\uinc$ is the {\em total wave} scattered by the obstacle under the sound-soft boundary condition $\utot|_\Gamma= 0$.

There are a number of ways to tackle problem \eqref{pde}. A popular option is to consider the problem on a bounded domain containing $\Gamma$ facilitated by an introduction of a transparent boundary condition. The exact transparent boundary condition is non-local in time and space and can be computed fast on special domains such as balls \cite{exact_special1,exact_special2,exact_special3,exact_special4}. \blue{Alternatively, an} approximate local boundary condition can be used, such as local absorbing boundary condition \cite{abc1,abc2},  methods based on the pole condition \cite{pole1,pole2} and perfectly matched layers \cite{pml}. All of these methods apply to convex domains of a special shape, e.g., circle or rectangle in 2D. This can result in an unnecessarily expensive method if $\Gamma$ is such that a very large convex domain is needed to encompass it, e.g., the elongated shape of an airplane. In such situations methods based on time-domain boundary integral equations are of advantage \cite{Costabel:2017,Sayas:2016,lbbook}. This particularly holds when the accurate computation of the far field potential is needed as the boundary integral potentials have the exact far-field behaviour encoded in their kernels.

Thus in this work we  represent the solution as a {\em single layer boundary integral potential}
\begin{equation}\label{eq:SLP}
    u(t,x) = S(\partial_t) \varphi(t,x) := \int_0^t\int_\Gamma k(t-\tau,|x-y|) \varphi(\tau,y) d\Gamma_y d\tau,
\end{equation}
where $\varphi \colon [0,T] \times \Gamma \to \R $ is an unknown density and $k$ is the fundamental solution which depends on the spatial dimension
\begin{equation}\label{eq:f_sol}
k(t,r) = \left\{
\begin{array}{cc}
   \displaystyle \frac{H(t-r)}{2\pi \sqrt{t^2-r^2}}  &  d = 2\vspace{.1cm}\\ 
    \displaystyle\frac{\delta(t-r)}{4\pi r} &  d = 3
\end{array}
\right.,
\end{equation}
where $H(\cdot)$ is the Heaviside function and $\delta(\cdot)$ the Dirac delta distribution. 

Taking the trace onto $\Gamma$ of \eqref{eq:SLP} we obtain {\em the boundary integral equation} for the unknown density $\varphi$: Find $\varphi$ such  that
\begin{equation}\label{eq:SLP_eq}
    V(\partial_t) \varphi(t,x) := \int_0^t\int_\Gamma k(t-\tau,|x-y|) \varphi(\tau,y) d\Gamma_y d\tau = g(t,x) \quad (t,x) \in [0,T] \times \Gamma.
\end{equation}
Once the density $\varphi$ is obtained, the solution $u$ can be recovered from \eqref{eq:SLP}.
All the numerical methods in this paper will be based on this boundary formulation of the scattering problem. 

Note that we have implicitly defined two operators above. The single layer potential $S(\partial_t) \cdot$ and single layer boundary integral operator $V(\partial_t) \cdot $· The motivation for the notation and the mapping properties are described in \cite{Lubich:1994,lbbook} and will be briefly explained in the next section.

A popular method for the discretization of \eqref{eq:SLP_eq} is the space-time Galerkin method originating in the work of Bamberger and Ha Duong \cite{BaHa:1986a,BaHa:1986b}. Alternatively, any spatial discretization can be combined with convolution quadrature (CQ) discretization in time as introduced by Lubich in \cite{Lubich:1994}. Combined with Galerkin discretization in space, the analysis of CQ discretization for standard problems of acoustics is nearly complete; see books \cite{Sayas:2016,lbbook}. \blue{Recently,} a combination with the method of fundamental solutions (MFS) was investigated in \cite{LH20_2859}. \blue{While} CQ discretization has a number of advantages -- powerful toolbox of analysis techniques and fast methods, relatively simple implementation, a vast literature on various applications, broadness of applicability -- a recognised weakness of the approach is its dissipative and dispersive properties \cite{Monk,lbbook}. While high-order methods based on Runge-Kutta can tackle this deficiency very effectively \cite{BaLuMe:2011,BaMeSc:2012}, it is still of interest to develop accurate (and cheaper) low-order variants of CQ. \blue{Despite the effectiveness of high-order methods based on Runge-Kutta time-stepping, low order methods are of interest due to their ease of implementation and ease of coupling with other solvers in bounded, inhomogeneous regions; see works on FEM-BEM coupling in the time-domain \cite{BaLuSa:2015,lbbook}. }  This has led to works on modified CQ schemes, see e.g., the works of Davies and Duncan \cite{DaDu:2013,DaDu:2014}, Weile \cite{Weile_mod}, and \cite{lb_ya_proc} by the second author. Here we make use of the method introduced by Weile. Our description will be closer to the description given in \cite{lb_torino} by the second author.  
The novelty of the current work is the description of the FFT based implementation of the scheme, the combination with MFS and numerous numerical experiments. \blue{The numerical results will show up to two orders of magnitude improvement in error when compared with the standard, low-order, schemes.} 
\blue{We should say that the modification we describe in this paper is applicable to high-er order schemes and are likely to improve accuracy even for high order Runge-Kutta schemes. The investigation of higher order schemes is however beyond the scope of the current work.}
We next give a brief content of the paper.

For the spatial discretization of \eqref{eq:SLP} and \eqref{eq:SLP_eq} we will consider two approaches: the Galerkin boundary element method (BEM) and the method of fundamental solutions (MFS). For the time discretization we will also consider two approaches: standard convolution quadrature (CQ) based on A-stable linear multistep methods and a class of modified CQ methods specially designed for wave propagation problems. We begin by describing standard CQ together with the modification mentioned above. Next, we describe spatial discretization approaches. Finally, we explain how the discretizations combine together and describe efficient algorithms for the implementation. The paper is finished by a section with extensive numerical experiments.

\section{One-sided convolutions and CQ}\label{sec:cq}

Convolution quadrature is a numerical method designed for the discretization of convolution
\[
K(\partial_t)g :=  \int_0^t k(t-\tau) g(\tau) d\tau,
\]
where $k$ is  a given kernel function with $K$ its Laplace  transform 
\[
K(s) =\{\mathscr{L} k\}(s) := \int_0^\infty e^{-st} k(t) dt \qquad \Re s > 0.
\]
If $k$ is not an integrable function but potentially a distribution with Laplace transform $K$  analytic for $\Re s > 0$ and \blue{bounded as}
\[
|K(s)| \leq C(\sigma) |s|^\mu \qquad \forall \Re s \geq \sigma > 0,
\]
\blue{for some $\mu \in \R $ and a positive function $C$.} 

we define the convolution $K(\partial_t) g$ via the inverse Laplace transform
\[
K(\partial_t)g(t) := \mathscr{L}^{-1} \{KG\}(t)  = \frac1{2\pi\mi} \int_{\sigma-\mi \infty}^{\sigma+\mi \infty} e^{st} K(s) G(s) ds,
\]
where $G = \{\mathscr{L} g\}$. As long as $G$ decays faster than $-\mu-1-\varepsilon$ for a fixed $\varepsilon > 0$, $K(\partial_t) g(t)$ is a continuous function that can be continuously extended to $t < 0$ by zero. Alternatively, if $g \in C^m(\mathbb{R}_{\geq 0})$ with $g(0) = g'(0) = \dots = g^{(m-1)}(0)= 0$ we have
\[
\left|\mathscr{L} \{g\}(s)\right| \leq |s|^{-m} \int_0^\blue{\infty}e^{-\sigma t} |g^{(m)}(t)| dt \qquad \Re s \geq \sigma > 0,
\]
and hence $K(\partial_t)g$ is continuous and $K(\partial_t)g(0) = 0$ if $m > \mu+1$. 

Given a time-step $\tstep > 0$ and a generating function of an A-stable linear multistep method\blue{, for instance}
\[
\delta(\zeta) = (1-\zeta)+\frac12(1-\zeta)^2 \;(\text{BDF2}) \qquad 
\delta(\zeta) = 2\frac{1-\zeta}{1+\zeta} \;\blue{(\text{trapezoidal rule}),}
\]
convolution quadrature (CQ) of the above convolution is defined by
\begin{equation}\label{eq:CQ_def}
K(\partial_t^\tstep)g(t) := \frac1{2\pi\mi} \int_{\sigma-\mi \infty}^{\sigma+\mi \infty} e^{st} K(s^\tstep) G(s) ds,
\end{equation}
where
\[
s^\tstep = \frac{\delta(e^{-s\tstep})}{\tstep}.
\]
The two linear multistep methods are second order and A-stable, which in terms of the generating functions is stated as
\begin{equation}\label{eq:LMSp}
\delta(e^{-z}) = z+\O(z^3)
\end{equation}
and
\begin{equation}
    \Re \delta(e^{-z})  > 0 \qquad \text{for } \Re z >0.
\end{equation}

For sufficiently smooth $g$ it is shown in \cite{Lubich:1994,FJStrap,lbbook} that 
\[
\left|K(\partial_t^\tstep)g(t)-K(\partial_t)g(t) \right| = \O(\tstep^2)
\]
for the above linear multistep methods. 

Let us note that $K(\partial_t^\tstep)g(t)$ is a discrete convolution, namely with $t_n = n \tstep$
\[
K(\partial_t^\tstep)g(t_n) = \sum_{j = 0}^n \omega_{n-j} g(t_j)
\]
where $\omega_j$ are the convolution weights defined by the generating function
\[
K\left(\frac{\delta(\zeta)}{\tstep}\right) = \sum_{j = 0}^\infty \omega_j \zeta^j.
\]
A crucial property for the stability of CQ is the composition rule which implies that solving the discrete equation
\[
K(\partial_t^\tstep)\varphi(t_n) = g(t_n) \qquad n = 0,\dots, N
\]
is equivalent to evaluating the discrete convolution
\[
\varphi_n = K^{-1}(\partial_t^\tstep)g(t_n) \qquad n = 0,\dots, N.
\]
\section{A modified CQ}\label{sec:mod_cq}

A fundamental property of solutions of the wave equation is the fact that they travel at finite speed. This is nicely visible in the definition of the fundamental solutions in \eqref{eq:f_sol} which are zero for $t < r$, i.e., are zero before the wave has arrived from the source to the receiver. This property is also visible in the Laplace domain since a simple calculation shows that
if $u$ is polynomially bounded and locally integrable then for $r \geq 0$
\begin{equation}\label{eq:fp_2freq}
u(t) \equiv 0 \text{ for } t < r  \; \implies \; 
\left|e^{sr}\mathscr{L}\{u\}(s)\right| \leq \int_0^\infty e^{-\sigma t} |u(t+r)| dt \text{ for } \Re  s \geq \sigma > 0.
\end{equation}
By the inverse Laplace transform and the Cauchy integral formula the reverse direction holds as well: For $r \geq 0$
\begin{equation}\label{eq:fp_2time}
    |K(s)e^{sr}| \leq C |s|^\mu \text{ for } \Re s > 0 \text{ and } \mu < -1 \implies
    k(t) \equiv 0 \text{ for } t \leq r.
\end{equation}

In particular if we know that the kernel  $K(s)$ is such that 
\[
\left|e^{st_m}K(s)\right| \leq C(\sigma) |s|^{\tilde \mu} \qquad \Re s \geq \sigma > 0,
\]
\blue{for some $\tilde \mu \in \R $ and a positive function $C$, then} $K(\partial_t)g(t) \equiv 0$ for $t \leq t_m$ for sufficiently smooth $g$. This motivates a modified CQ, where the quadrature is only applied after $t > t_m$:
\[
K(\partial_t^\tstep;t_m) g(t) := \left\{
\begin{array}{cc}
\displaystyle\frac1{2\pi\mi} \int_{\sigma-\mi \infty}^{\sigma+\mi \infty} e^{s(t-t_m)} e^{s^\tstep t}K(s^\tstep) G(s) ds,     & t > t_m  \\
0     &    t \leq t_m.
\end{array}
\right.
\]
This is again a discrete  convolution
\begin{equation}\label{eq:mod_cq}
K(\partial_t^\tstep;t_m) g(t_n)= \sum_{j = 0}^n \omega_{n-j; m} g(t_j)
\end{equation}

where the generating function for the weights is given by
\begin{equation}\label{eq:mod_cq_weights}
\zeta^m e^{t_m\frac{\delta(\zeta)}{\tstep}} K\left(\frac{\delta(\zeta)}{\tstep}\right) =\zeta^m e^{m\delta(\zeta)} K\left(\frac{\delta(\zeta)}{\tstep}\right) = \sum_{j = 0}^\infty \omega_{j;m} \zeta^j.
\end{equation}
Note that $\omega_{j;m} = 0$ for $j < m$ and hence $K(\partial_t^\tstep;t_m) g(t) = 0$ for $t < t_m$. This modified CQ avoids approximating values that we know are zero and furthermore, as we will see later, it has the potential to give much more accurate solutions in certain cases.

\section{Spatial discretization of integral operators}

Time-domain boundary integral operators fit seamlessly into the setting of Section~\ref{sec:cq}.
Namely, the Laplace domain single layer potential is given by
\[
S(s) \varphi(x) := \int_\Gamma K(s,|x-y|) \varphi(y) d\Gamma_y \qquad x \in \R^d \setminus \Gamma,
\]
where $K(s,r)$ is the Laplace transform of the free space Green's function \blue{for} wave equation in \eqref{eq:f_sol}
\[
K(s,r) =  \left\{
\begin{array}{cc}
   \displaystyle \frac1{2\pi}K_0(sr)  &  d = 2\vspace{.1cm}\\ 
    \displaystyle\frac{e^{-sr}}{4\pi r} &  d = 3
\end{array}
\right.,
\]
with $K_0(\cdot)$ the modified Bessel function of the second kind.
From the definition of $K$ for $d = 3$ and from the asymptotic behaviour of $K_0$ \cite{NIST} we have that for $r > 0$ there exists a constant $C(r) > 0$ such that
\begin{equation}\label{eq:esrK_bound}
    |e^{sr}K(s,r)| \leq C(r) \qquad \text{for } \Re s > 0.
\end{equation}
Recall that this corresponds to the time-domain property $k(t,r) = 0$ for $t < r$; see \eqref{eq:f_sol}.

\blue{The single} layer potential is continuous across the boundary, hence in the frequency domain we are interested in the single-layer boundary integral operator $V(s) \colon H^{-1/2}(\Gamma) \to  H^{1/2}(\Gamma)$
\[
V(s) \varphi(x) := \int_\Gamma K(s,|x-y|) \varphi(y) d\Gamma_y \qquad x \in \Gamma.
\]
By now, \cite{BaHa:1986a,Sayas:2016,lbbook}, it is well-known that $V(s)$ and $V^{-1}(s)$ are bounded for $\Re s > 0$ as 
\[
\begin{split}
    \|V(s)\|_{H^{-1/2}(\Gamma) \to H^{1/2}(\Gamma)} &\leq C \max(1,|s|^{-2}) \frac{|s|}{\Re s}\\
    \|V^{-1}(s)\|_{H^{1/2}(\Gamma) \to H^{-1/2}(\Gamma)} &\leq C \max(1,|s|^{-1}) \frac{|s|^2}{\Re s}.
\end{split}
\]

We now describe two standard approaches to discretizing these operators.

\subsection{Galerkin boundary element methods (BEM)}\label{sec:Gal}

Let $\Gamma$ be subdivided into $M$ non-overlapping panels $\Gamma_i$, $i = 1, \dots, M$. We denote by $X_h \subset H^{-1/2}(\Gamma)$ the finite-dimensional space of piecewise constant functions. A basis of the space is given by $\varphi_1, \dots, \varphi_M$ defined as
\[
\varphi_i(x) = \left\{\begin{array}{cc}
    1 & x \in \Gamma_i \\
    0 & \text{otherwise}.
\end{array}\right.
\]
Then we define the matrix representation of the Galerkin discretization of the single-layer operator by
\begin{equation}
\begin{split}
    \left(\Vg(s)\right)_{ij} &:= \int_\Gamma\int_\Gamma K(s,|x-y|) \varphi_i(y) \varphi_j(x) d\Gamma_y d\Gamma_x\\
    &= \int_{\Gamma_i} \int_{\Gamma_j}K(s,|x-y|) d\Gamma_y d\Gamma_x.
\end{split}
\end{equation}
An important property of this matrix, implied by \eqref{eq:esrK_bound}, is that 
\[
\left|e^{s r_{ij}}\left(\Vg(s)\right)_{ij} \right| \leq const
\]
for $\Re s \geq \sigma \geq 0$, $r_{ij} > 0$, and 
\[
r_{ij} \leq \dist(\Gamma_i,\Gamma_j).
\]
Again, we can interpret this in the time-domain as the fact that information travelling from one panel, e.g., $\Gamma_i$, needs time greater than $r_{ij}$ to reach the second panel, $\Gamma_j$.

With this spatial discretization in place, the semi-discretization of \eqref{eq:SLP_eq} reads: Find $\boldsymbol{\varphi}(t) \in \mathbb{R}^M$ such that
\begin{equation}\label{eq:Gal_semi}
    \Vg(\partial_t) \boldsymbol{\varphi} (t) = \mathbf{g}(t),
\end{equation}
where the right-hand side is the projection $\mathbf{g}(t) \in \mathbb{R}^M$ of the boundary data defined by
\[
\left(\mathbf{g}(t)\right)_i = \int_\Gamma g(t,x) \varphi_i(x) d\Gamma_x.
\]

\subsection{Method of Fundamental Solutions}

As we have seen in the previous sections, the Laplace domain solution $U(s) = \mathscr{L}\{ u\}(s)$ can be represented by the single layer potential:
\[
U(s,x) = \int_\Gamma K(s,|x-y|) \Phi(s,y) d\Gamma_y,
\]
where we have now explicitly given the dependence of the density $\Phi$ on the Laplace domain parameter $s$. \blue{The }Galerkin discretization described above, while accurate, stable and well-understood is expensive and non-trivial to implement due to the singular, double integrals. A simple alternative is the Method of Fundamental Solutions (MFS) \cite{BarB,FairK}. In the MFS, we approximate $U(s)$ by a sum of {\em sources} located at $y_j$, $j = 1,\dots, K$ 
\begin{equation}
    U(s,x) \approx \sum_{j = 1}^K \Phi_j(s) K(s,|x-y_j|).
\end{equation}
The unknown coefficients are determined by solving the least squares problem derived from the boundary condition: Find $\Phi_j(s)$ such that
\begin{equation}
     \sum_{j = 1}^K \Phi_j(s) K(s,|x_i-y_j|)  = G(s,x_i), \quad i = 1,\dots,M,
\end{equation}
where $x_i \in \Gamma$ are the {\em collocation} points and $G = \mathscr{L}\{g\}$. Note that the source points $y_j$ are located in the exterior of the computational domain and certainly away from the boundary, hence no singularities occur in the system.
A careful numerical analysis of this method for
 connected planar \blue{domains} with analytic boundary has been given \blue{in }\cite{BarB}.

Again, we can define a system matrix $\Vmfs(s) \in \mathbb{C}^{M \times K}$
\begin{equation}\label{eq:V_MFS}
    \left(\Vmfs(s)\right)_{ij} = K(s,|x_i-y_j|),
\end{equation}
which in contrast with the Galerkin method is rectangular and trivial to compute. Importantly again we have the property
\[
\left|e^{s r_{ij}}\left(\Vmfs(s)\right)_{ij} \right| \leq const,
\]
where now $r_{ij} = |x_i-y_j|$.

The semi-discrete system again reads: Find $\boldsymbol{\varphi}(t) \in \mathbb{R}^K$ such that
\begin{equation}\label{eq:MFS_semi}
    \Vmfs(\partial_t) \boldsymbol{\varphi} (t) = \mathbf{g}(t),
\end{equation}
where the right-hand side is now simply
\[
\left(\mathbf{g}(t)\right)_i = g(t,x_i), \qquad i = 1,\dots,M.
\]
\section{Accuracy of the two time-discretizations}
The approximation due to the standard CQ, according to \eqref{eq:CQ_def} consists of replacing in the Laplace domain the kernel $K(s,r)$ by its approximation $K(s^{\tstep},r)$ with $s^{\tstep} = \delta(e^{-s\tstep})/\tstep$. In particular, in 3D this means that we are replacing $e^{-sr}$ by $e^{-s^{\tstep}r}$. The error of this approximation can be bounded as
\[
\left|e^{-sr}-e^{-s^{\tstep}r}\right| \leq Cr|s-s^{\tstep}| \leq C|rs||s\tstep|^2,
\]
for some constant $C> 0$ that is allowed to change from one step to another and $|s\tstep| < c$ for some small enough constant $c >0$. \blue{Since $s = i \omega,$ if} we denote by $\omega_{\text{max}}$ the largest frequency present in the system, \blue{the estimate above implies that} to get a fixed accuracy we expect that we \blue{need
$r\tstep^2\omega_{\text{max}}^3$ to be} small enough, i.e., the time-step $\tstep$ needs to be chosen proportional to $\omega_{\text{max}}^{-3/2}$ resulting \blue{in the number} of degrees of freedom per wavelength not being fixed.

Instead, in the modified method the approximation is of the form $e^{-st_m}e^{-s^{\tstep}(r-t_m)}$, where $t_m$ is the largest time-step smaller than $r_{ij} \leq r$. The approximation is now
\[
\left|e^{-sr}-e^{-st_m}e^{-s^{\tstep}(r-t_m)}\right|
\leq \left|e^{-s(r-t_m)}-e^{-s^{\tstep}(r-t_m)}\right|
\leq C|(r-t_m)s||s\tstep|^2.
\]
In both the case of the Galerkin and MFS discretizations, $t_m$ can be chosen so that
$r-t_m = \O(h)+\O(\tstep)$, where $h$ is the spatial mesh-width. If, as is usually the case, we have $h \propto \tstep$, we obtain that the above error remains at a fixed accuracy if 
$\tstep^3\omega_{\text{max}}^3$ is small enough, i.e., a fixed number of degrees of freedom is sufficient.

Similar arguments hold for the 2D case, recalling the asymptotic behaviour of $K_0$ for large arguments \cite[10.25.3]{NIST}:
\[
K_0(z) \sim \sqrt{\frac{\pi}{2z}} e^{-z} \qquad \text{as } z \to \infty, \; |\arg z | < \pi.
\]

\section{Fully discrete system and algorithmic implementation}
In the \blue{following,} $\Vb(s)$ will denote either the Galerkin system matrix $\Vg$ or the corresponding matrix $\Vmfs$ obtained via MFS. Hence, the semi-discrete system reads: Find $\boldsymbol{\varphi}(t)$ such that 
\[
    \Vb(\partial_t) \varphib (t) = \mathbf{g}(t).
\]

\blue{The fully} discrete system can directly be obtained by standard convolution quadrature discretization in time: Find $\varphib^{\tstep}_j$, $j = 0,\dots, N$, such that
\begin{equation}\label{eq:fully_CQ}
    \Vb(\partial_t^\tstep) \varphib^{\tstep} (t_n) = \mathbf{g}(t_n), \qquad
    n = 0,\dots,N,
\end{equation}
where $t_j = j\tstep$, $\tstep > 0$ is the time-step and $t_N = T$ the final time.

To describe the modified scheme, recall that we are given distances $r_{ij}$ such that
\[
\left|e^{s r_{ij}}\left(\Vb(s)\right)_{ij} \right| \leq const.
\]
Next we define
\[
m_{ij} = \left\lfloor \frac{r_{ij}}{\tstep}\right\rfloor
\]
so that $0 \leq t_{m_{ij}} \leq r_{ij}$ and is the largest time-step satisfying this inequality.

The modified scheme then reads: Find $\widetilde{\varphib}^{\tstep}_\ell$, $\ell = 0,\dots, N$, such that
\begin{equation}\label{eq:fully_mod}
    \Vb(\widetilde{\partial}^\tstep_t) \widetilde{\varphib}^{\tstep} (t_n) = \mathbf{g}(t_n), \qquad
    n = 0,\dots,N,
\end{equation}
where 
\begin{equation}\label{eq:fully_mod_conv}
    \Vb(\widetilde{\partial}^\tstep_t) \widetilde{\varphib}^{\tstep} (t_n)
:= \sum_{\ell = 0}^n \widetilde{\omega}_{n-\ell} \widetilde{\varphib}^{\tstep}_\ell
\end{equation}
and the modified weights are given by the following generating function, c.f.\ \eqref{eq:mod_cq_weights},
\begin{subequations}\label{eq:fully_mod_def}
\begin{align}
\widetilde{\Vb}(\zeta)    &=
\sum_{n = 0}^\infty \widetilde{\omega}_{n} \zeta^n.\label{eq:fully_mod_defa}\\
\left(\widetilde{\Vb}(\zeta)\right)_{ij} &:=\zeta^{m_{ij}} e^{m_{ij}\delta(\zeta)} \left(\Vb\left(\frac{\delta(\zeta)}{\tstep}\right)\right)_{ij}.
\label{eq:fully_mod_defb}
\end{align}
\end{subequations}
If $\tilde \omega_0$ is invertible in the least squares sense, then the solution is given by the iteration: Find $\widetilde{\varphib}^{\tstep}_n$ such that
\begin{equation}\label{eq:fully_mod_MOT}
\widetilde{\omega}_{0}  \widetilde{\varphib}^{\tstep}_n
= \gb(t_n)-\sum_{\ell = 0}^{n-1} \widetilde{\omega}_{n-\ell} \widetilde{\varphib}^{\tstep}_\ell
\end{equation}
in the least squares sense.

Note that the weights in \eqref{eq:fully_mod_conv} are defined by \blue{the} generating function \eqref{eq:fully_mod_defa} as in the standard CQ. Hence, all the FFT based algorithms designed for CQ are still available. For completeness, we \blue{describe the} most basic of these algorithms. More \blue{details} and other algorithms can be found in \cite{lbbook,BaSc:2012}.

First of \blue{all,} notice that \eqref{eq:fully_mod_defa} is a Taylor expansion of the function $\widetilde{\Vb}(\zeta)$ analytic for $|\zeta| < 1$. \blue{Hence,} we can represent the coefficients by the Cauchy integral formula
\[
\begin{split}
\widetilde{\omega}_j &= \frac1{2\pi\mi} \oint_{\mathcal{O}} \widetilde{\Vb}\left(\zeta\right)\zeta^{-j-1} d\zeta\\
&= \lambda^{-j} \int_0^1 \widetilde{\Vb}\left(\lambda e^{2\pi \mi \theta}\right)e^{-2\pi \mi j \theta} d\theta,
\end{split}
\]
where \blue{we have chosen the disk of radius $0 <\lambda < 1$ as the contour}. Applying the composite trapezoidal rule to this (periodic) integral we obtain an approximation of the weights
\begin{equation}\label{eq:weights_app}
    \widetilde{\omega}_j \approx \frac{\lambda^{-j}}{N+1} \sum_{\ell = 0}^N
    \widetilde{\Vb}\left(\lambda \zeta_{N+1}^{-\ell}\right)\zeta_{N+1}^{\ell j},
\end{equation}
where $\zeta_{N+1} = e^{\frac{2\pi\mi}{N+1}}$. This approximation is valid for $j = 0,\dots,N$ and the error is of $\O(\lambda^{N+1})$. In finite precision \blue{arithmetic,} due to the multiplication with the factor $\lambda^{-j}$, $\lambda$ cannot be chosen too small. In practice $\lambda = \varepsilon^{1/2(N+1)}$, where $\varepsilon$ is close to machine precision, is a good choice giving as error $\sqrt{\varepsilon}$; for details see \cite{Lubich:1994} and for how to improve on this accuracy see \cite{lbbook}. 

Importantly, in \eqref{eq:weights_app} we recognise the inverse discrete Fourier transform, which can be computed in $\O(N\log N)$ time for all $j = 0,\dots,N$ by the fast Fourier transform (FFT). In a similar \blue{way,} further FFT based algorithms as described in \cite{lbbook} are available for the modified scheme. Note, however, that data sparse techniques such as FMM \cite{fmm:93,Greengard_phd} and $\mathcal{H}$-matrices \cite{Hackbook} cannot be directly applied to the computation of the spatial operator $\Vb(s)$ due to the discontinuous changes in the spatial kernel resulting from \eqref{eq:fully_mod_defb}. \blue{An alternative would be to use the exact shifts for {\em groups} of panels and their interaction, rather than for each pair in the Galerkin matrix. In this way, the kernel would still be smooth allowing for data sparse techniques and the new method would still significantly improve the accuracy. Detailed research would be needed to determine the optimal way of choosing these groups of panels. This is similar to the approach taken in \cite{BaKa:2014a}, where the near field was treated in the time-domain and the far field was treated in the Fourier domain with data sparse techniques. Another possibility would be to smooth out the shifts similar to what was done in \cite{SaVe:2013} in the case of the space-time Galerkin method.}

In the numerical experiments of this paper, we will use a parallel FFT algorithm introduced in \cite{BaSa:2008}, whose details we describe next. First of all we note that the approximation 
\eqref{eq:weights_app} of the weights $\widetilde\omega_j$ is valid also for $j = -N,\dots,-1$ if we define $\widetilde\omega_k = 0$ for $k < 0$. Furthermore, extending the sum in \eqref{eq:fully_mod_conv} to $N$ we have that the system to be solved is
\[
\sum_{\ell = 0}^N \widetilde{\omega}_{n-\ell} \widetilde{\varphib}^{\tstep}_\ell = \mathbf{g}_n
\]
for $n = 0,\dots,N$. Substituting now the approximation \eqref{eq:weights_app}, we have the following (approximate) system to be solved
\[
\sum_{\ell = 0}^N \left[\frac{\lambda^{-(n-\ell)}}{N+1} \sum_{k = 0}^N
    \widetilde{\Vb}\left(\lambda \zeta_{N+1}^{-k}\right)\zeta_{N+1}^{k (n-\ell)}\right] \widetilde{\varphib}^{\tstep}_\ell \approx \mathbf{g}_n.
\]
Rearranging gives
\[
\frac{\lambda^{-n}}{N+1}\sum_{k = 0}^N \widetilde{\Vb}\left(\lambda \zeta_{N+1}^{-k}\right)\left[ \sum_{\ell = 0}^N
     \lambda^{\ell}\widetilde{\varphib}^{\tstep}_\ell \zeta_{N+1}^{-\ell k }\right]\zeta_{N+1}^{k n }\approx \mathbf{g}_n.
\]
Multiplying both sides by $\lambda^n$ and applying the discrete Fourier transform we obtain $N+1$ decoupled linear systems to be solved:
\begin{equation}\label{eq:decoupled}
\widetilde{\Vb}\left(\lambda \zeta_{N+1}^{-k}\right)\left[ \sum_{\ell = 0}^N
    \lambda^{\ell}\widetilde{\varphib}^{\tstep}_\ell \zeta_{N+1}^{-\ell k } \right]\approx  \sum_{n = 0}^{N}\lambda^n\mathbf{g}_n \zeta_{N+1}^{-kn},  
\end{equation}
$k = 0,\dots,N$. Thus, once the (scaled) discrete Fourier transform of the data $\mathbf{g}_n$ is computed, $N+1$ linear systems need to be solved in parallel. The solution $\varphi_\ell^{\tstep}$ is then recovered by another application of the inverse discrete Fourier transform. Note that using symmetry, roughly half of the problems in \eqref{eq:decoupled} are conjugates of the other half and hence need not to be solved.

\begin{remark}
Note that if $\tstep < \min_{i,j} r_{i,j}$ then $m_{ij} > 0$ for all $i,j$ and  from the definition of the weights $\tilde \omega_n$ in \eqref{eq:fully_mod_def} it follows that 
\[
\left(\tilde \omega_n \right)_{ij} = \left(\widetilde{\Vb}(0)\right)_{ij} 
= 0
\]
for all $i,j$. In this case we are not able to solve \eqref{eq:fully_mod} by using the iteration \eqref{eq:fully_mod_MOT}. However, as we will see in the numerical experiments, the formula \eqref{eq:decoupled} is feasible even in this case and gives good results. 
\end{remark}

\section{Numerical experiments}

\subsection{Transient wave  scattering and MFS}

We first consider exterior scattering by either the unit disk $\mathcal{D} = \{ x\ \in \mathbb{R}^2 \;:\; |x| < 1\}$ or the union of two open ellipses $\mathcal{E} = \mathcal{E}_1 \cup \mathcal{E}_2$, where the ellipses are defined using the identification \blue{of $\mathbb{R}^2$} with the complex plane via
\[
\mathcal{E}_1 := \left\{f(z)-2 \;:\; |z| < 1\right\}\qquad
\mathcal{E}_2 := \left\{f(z)+2 \;:\; |z| < 1\right\},
\]
with $f$ the conformal map
\[
f(z) =  \frac12 e^{\mi \frac{\pi}{2}}(z+\frac15 z^{-1}).
\]
As the Dirichlet data in both cases we use
\[
g(\red{t,x}) = \sin(\omega(t-x \cdot \alpha)) \red{g}(t-4-\blue{x \cdot \alpha})
\]
with
\[
g(t) = e^{-(t/0.7)^2}, \qquad \alpha = (0, -1)^T
\]
and $\omega$ a parameter that we will choose as either $\omega= 1$ or $\omega = 5$.

For the Galerkin boundary element method, we let $X_h$ be the space of piecewise constant boundary element functions. Whereas for the method of fundamental solutions in the case of the unit disk we let the collocation points be
\[
x^{\mathcal{D}}_i = e^{2\pi \mi i/M}, \qquad i = 0,\dots, M-1  \qquad (\mathcal{D})
\]
and the source points
\[
y_j^{\mathcal{D}} = R e^{2\pi \mi j/K}, \qquad j = 0,\dots, K-1  \qquad (\mathcal{D})
\]
for some $R < 1$. For all exterior scattering problems we use $R = 0.9$. \red{At the end of the section, we will also show results for the interior problem for the unit disk where we choose $R=1.1$.}

\begin{figure}
    \centering
    \includegraphics[width=.4\textwidth]{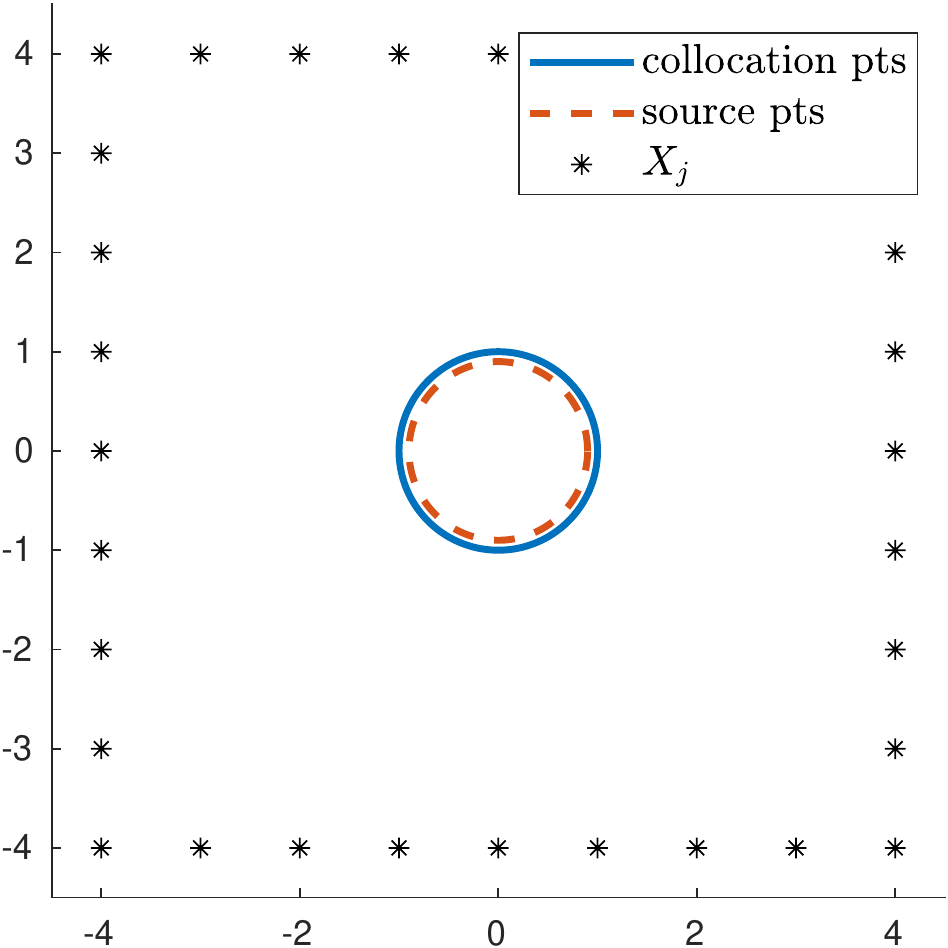}
    \hspace{.25cm}
    \includegraphics[width=.4\textwidth]{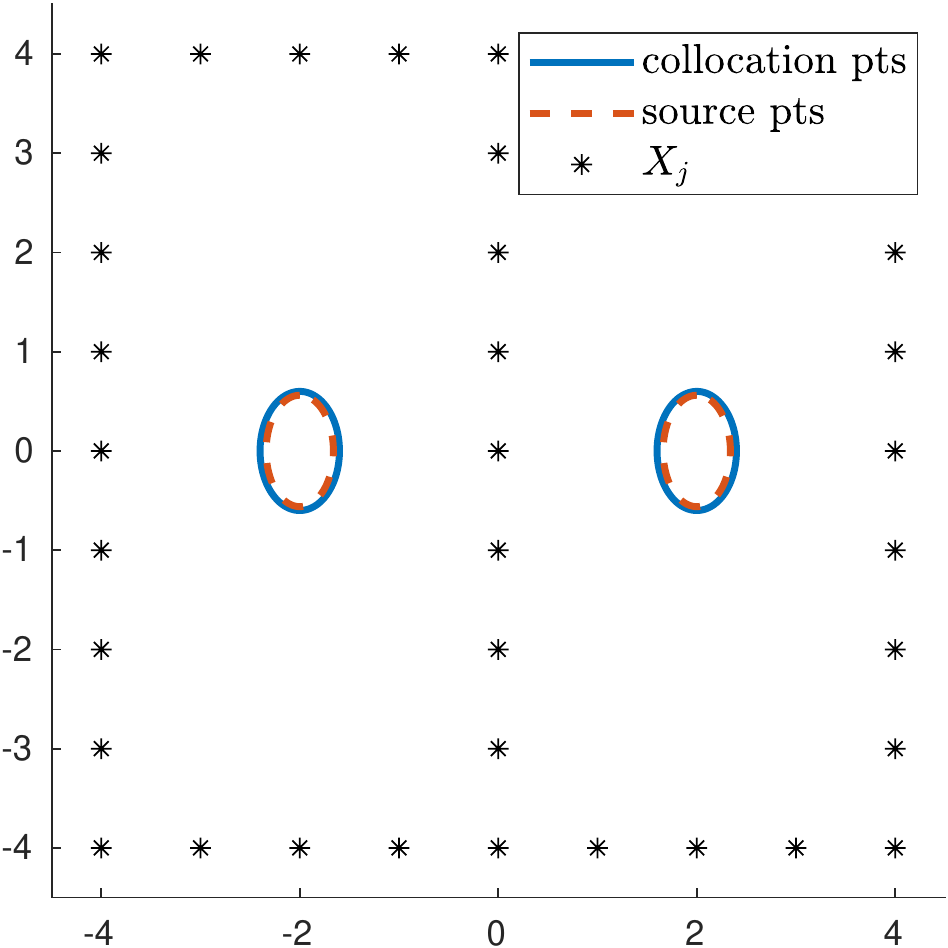}
    \caption{We show the position of the collocation points (solid lines), position of the source points (dashed lines) and the points $X_j$ where we evaluate the error. On the left the setting for the exterior scattering by the unit disk is shown, whereas on the right by two ellipses. }
    \label{fig:domains_ext}
\end{figure}

For each of the two ellipses we make use of the conformal map $f$ and define the collocation points to be
\[
x^{\mathcal{E}_1}_i = f(x^{\mathcal{D}}_i)-2, \qquad i = 0,\dots, M-1  \qquad (\mathcal{E}_1)
\]
and source points
\[
y^{\mathcal{E}_1}_j = f(y^{\mathcal{D}}_j)-2, \qquad j = 0,\dots, M-1  \qquad (\mathcal{E}_1)
\]
and analogously for $\mathcal{E}_2$. For a plot of the source and collocation points in each of the two cases see Figure~\ref{fig:domains_ext}.

This explains the spatial discretization of the problems. In time, we either use the standard CQ based either on BDF2 or the trapezoidal rule or its modified version. The fully discrete system is solved using the formula \eqref{eq:decoupled}, where each linear system is solved as a least squares problem.

Finally, to measure the error we choose some test points $X_\ell$, $\ell = 1,\dots,L$ placed in the exterior domain. As we do not have an exact solution at hand, we will use a finer mesh to obtain an accurate approximation denoted by $u_{\text{ex}}$. The error measure is \blue{the} maximum error over all test points and all time-steps
\begin{equation}\label{eq:error}
\text{error} = \max_{\ell,n} |u^h(\red{t_n,X_\ell})-u_{\text{ex}}(\red{t_n,X_\ell})|,
\end{equation}
where $u^h$ denotes the numerical solution with a coarser mesh in time and space.

\begin{figure}
    \centering
    \includegraphics[width=.8\textwidth]{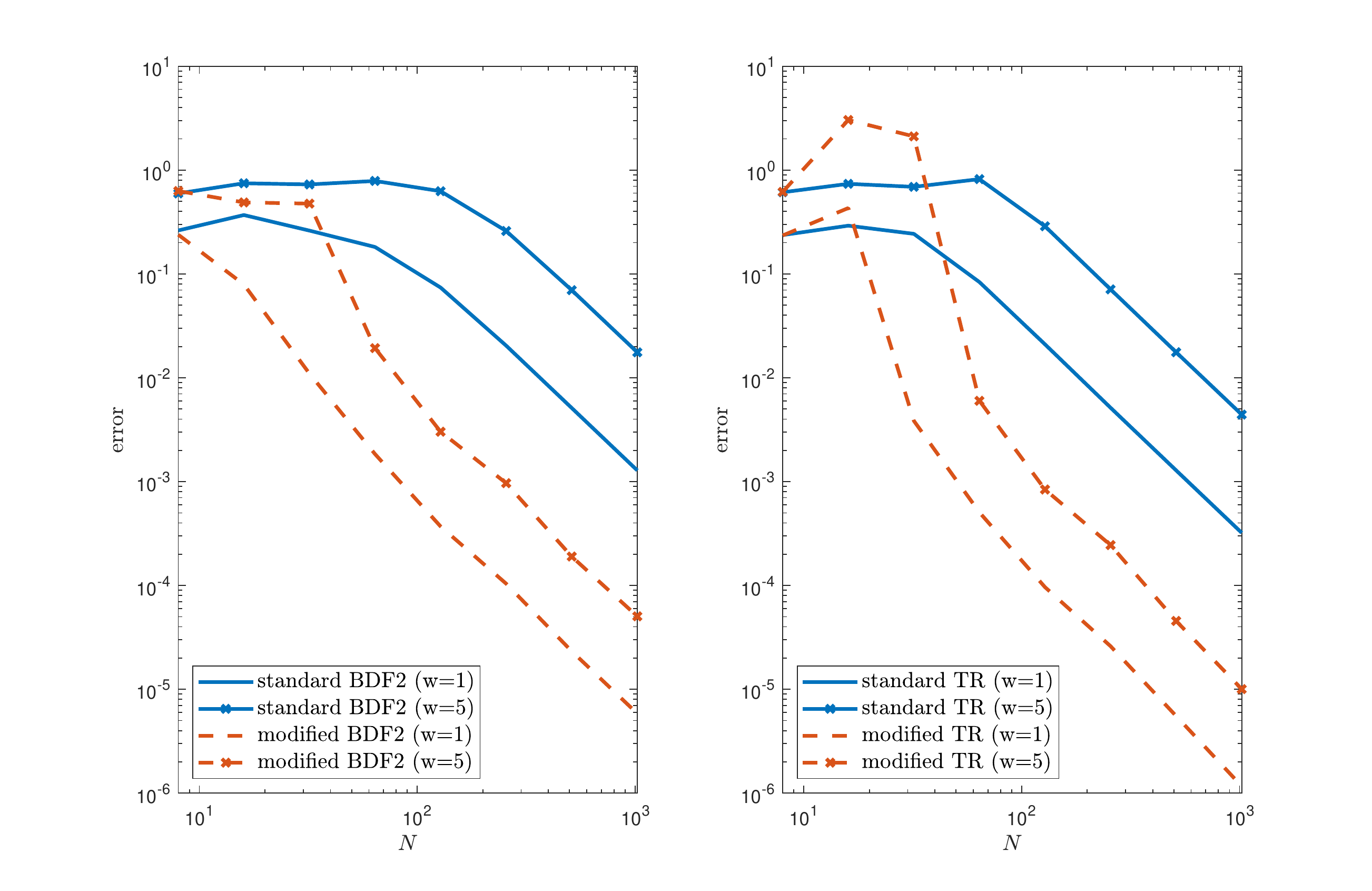}
    \caption{Convergence of the MFS with BDF2 based methods on the left and MFS with trapezoidal rule based methods on the right for the case of exterior scattering by the unit disk. }
    \label{fig:conv_disk_ext}
\end{figure}

In Figure~\ref{fig:conv_disk_ext}, in the case of the exterior scattering by the unit disk, we plot the convergence of the error for the standard CQ based on BDF2 and trapezoidal rule and its modified counterparts.
In these calculations we used $M = 2000$ collocation points and $K = 1000$ source points for the MFS to compute $u^h$. To compute the accurate solution $u^{\text{ex}}$ we used $M^{\text{ex}} = 3000$ collocation points and $K^{\text{ex}} = 1500$ source points and in time the modified CQ based on BDF2 with $N^{\text{ex}} = 2^{12}$ time-steps. We see a significant improvement when using the modified method.

\begin{figure}
    \centering
    \includegraphics[width=.8\textwidth]{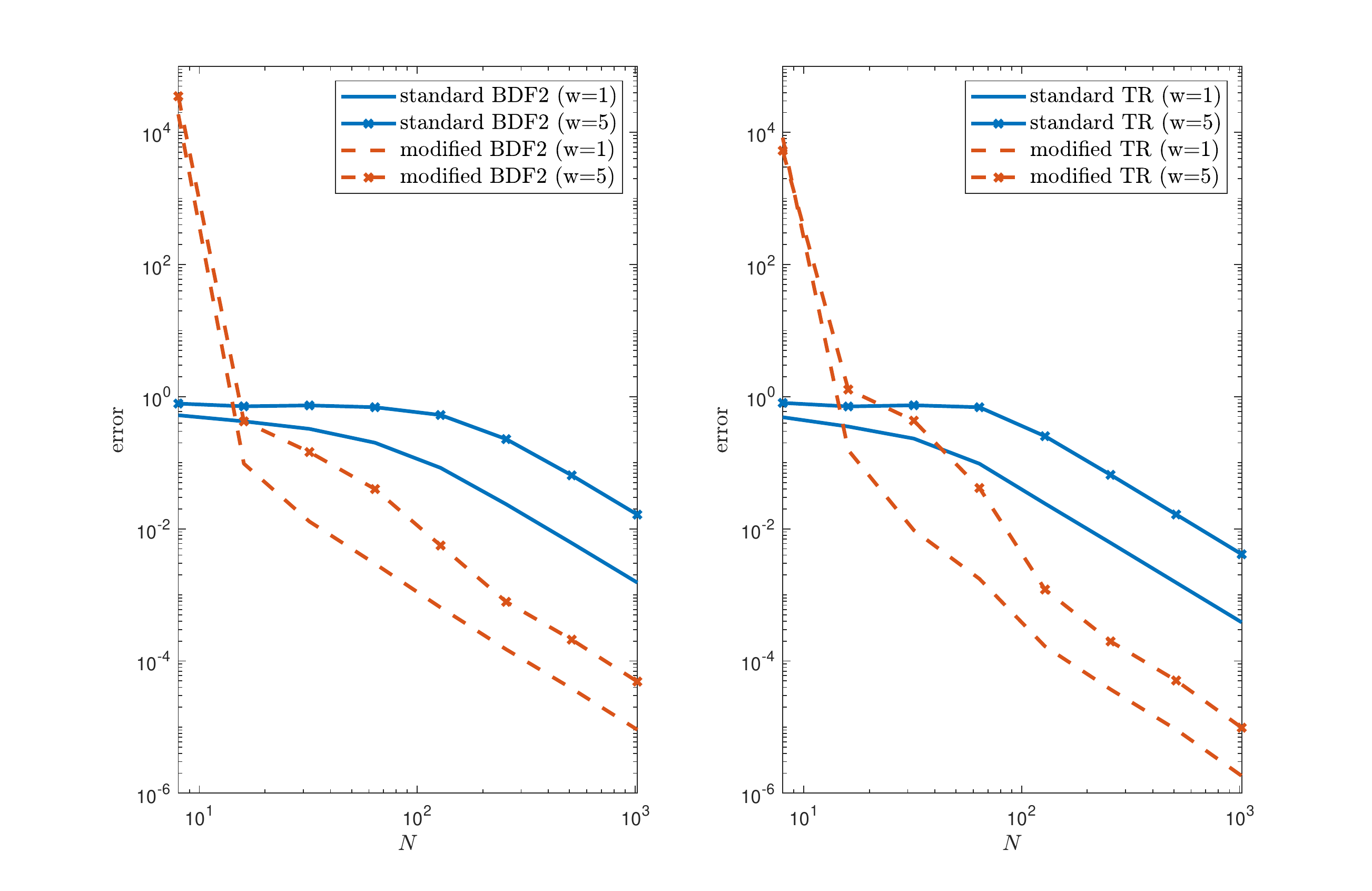}
    \caption{Convergence of the MFS with BDF2 based methods on the left and MFS with trapezoidal rule based methods on the right for the case of exterior scattering by two ellipses. }
    \label{fig:conv_ellipses}
\end{figure}

For the exterior scattering by two ellipses we plot the convergence of the error for the standard CQ based on BDF2 and trapezoidal rule and its modified counterparts in Figure~\ref{fig:conv_ellipses}. In these calculations we used $M = 4000$ collocation points and $K = 2000$ source points for the MFS to compute $u^h$. To compute the accurate solution $u^{\text{ex}}$ we used $M^{\text{ex}} = 6000$ collocation points and $K^{\text{ex}} = 3000$ source points and in time the modified CQ based on BDF2 with $N = 2^{12}$ time-steps. We see a significant improvement when using the modified method. \blue{Note, however, that for the modified scheme there is a transient region where the error is very large; see Figure~\ref{fig:conv_ellipses}. We believe, though we have no proof,  that the reason for this is that the new scheme does not damp high frequencies whereas CQ does. This makes CQ stable even when it is inaccurate, whereas the modified scheme produces controlled results only once it starts converging.}

%============= ellipses figures ===============
\begin{figure}[htb]%\hspace{-25mm}
\centering
  \begin{subfigure}[b]{0.33\textwidth}
    \centering
    \includegraphics[width=\linewidth]{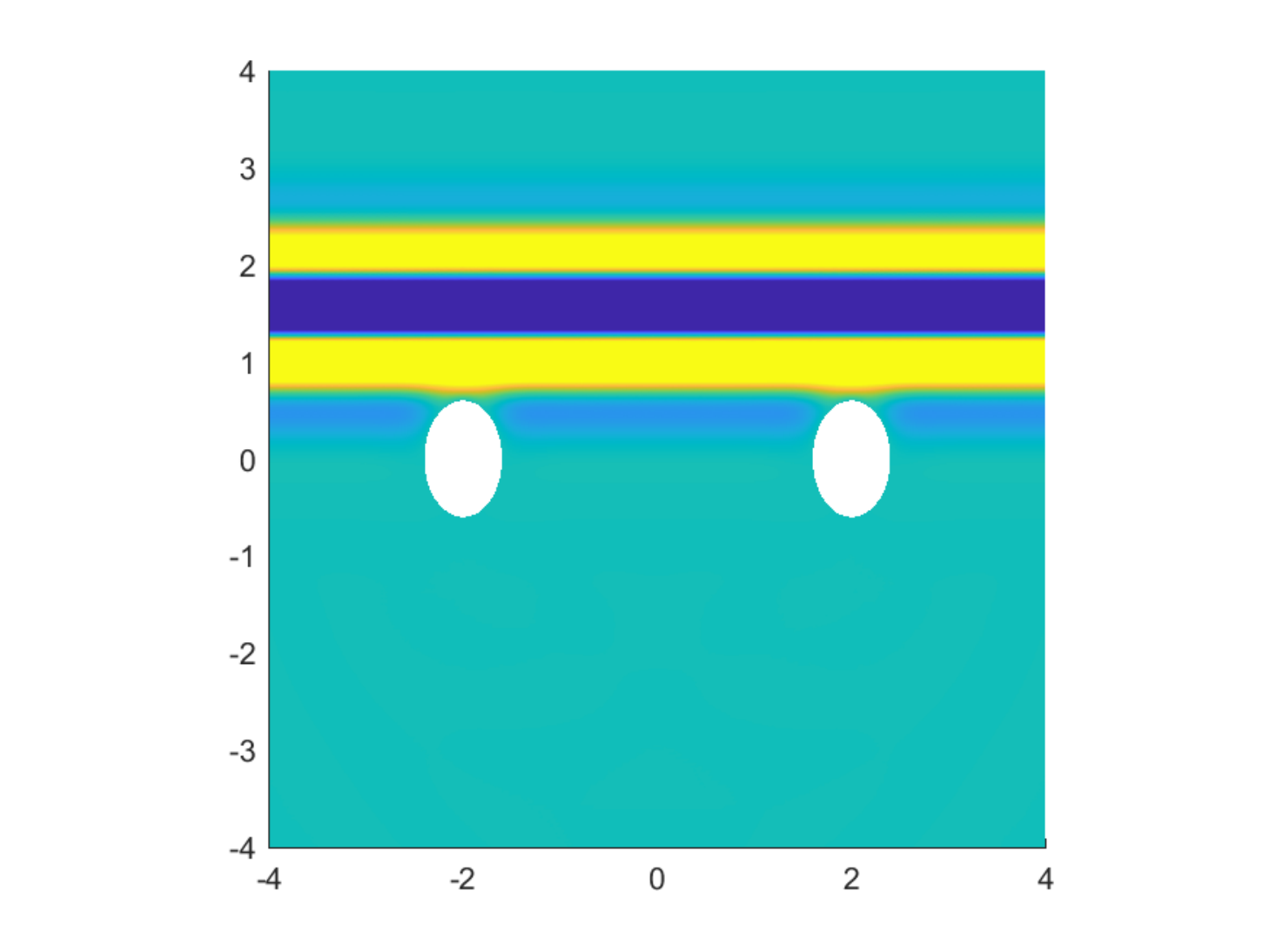}
    \caption{At t=2.5}
  \end{subfigure}\hfil%   
  \begin{subfigure}[b]{0.33\textwidth}
    \centering
    \includegraphics[width=\linewidth]{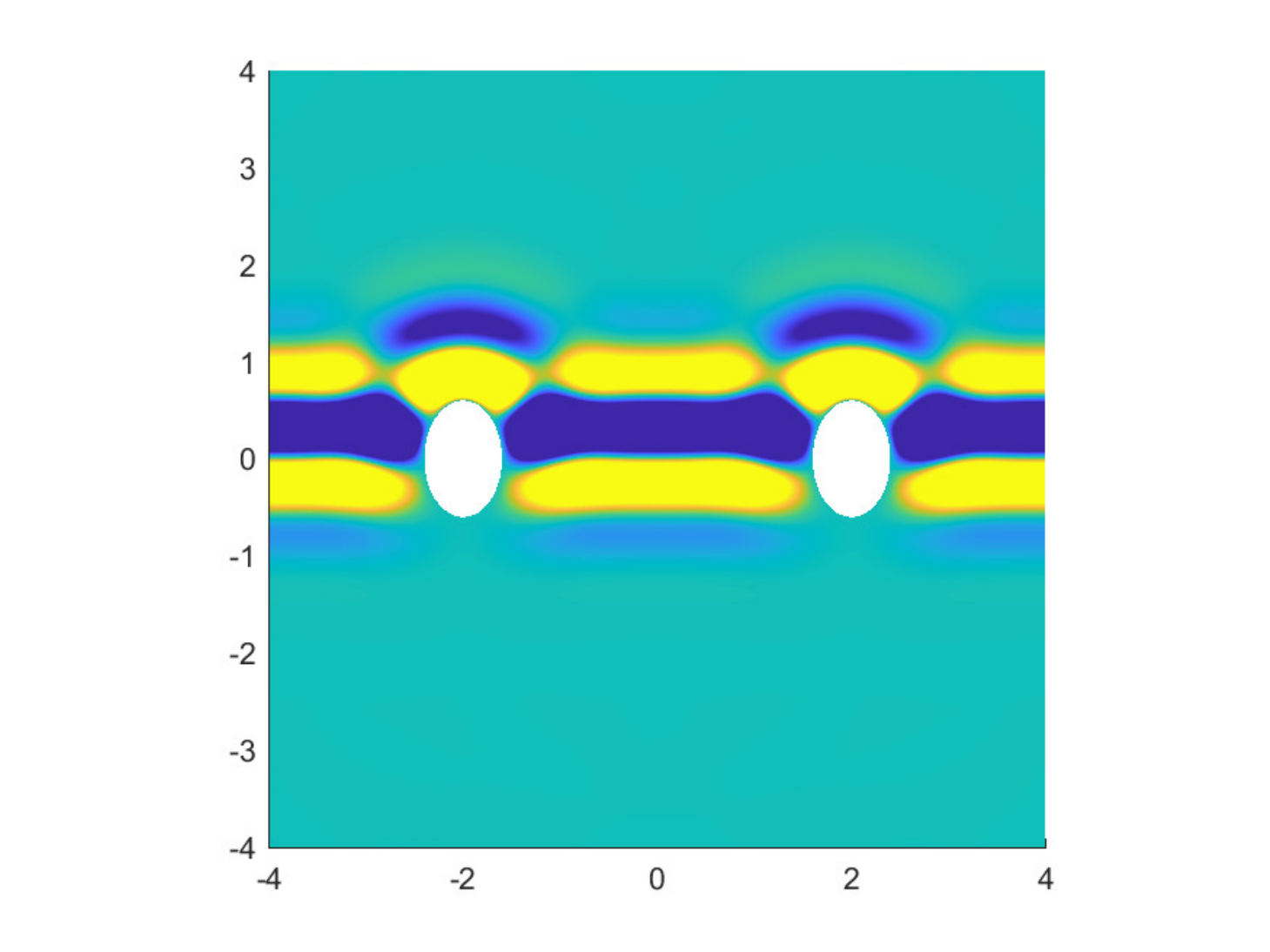}
    \caption{At t=3.75}
  \end{subfigure}\hfil%  
  \begin{subfigure}[b]{0.33\textwidth}
    \centering
    \includegraphics[width=\linewidth]{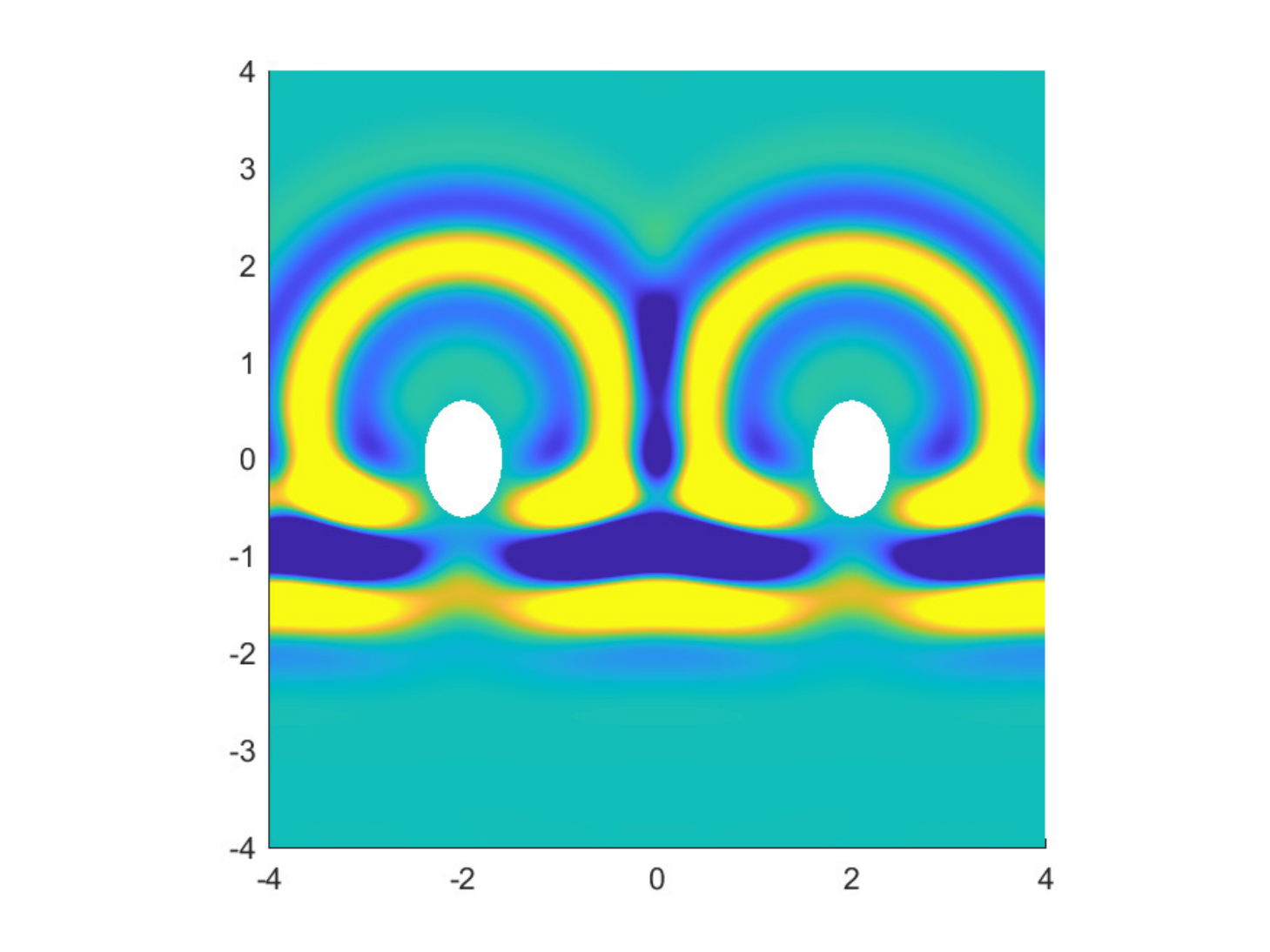}
    \caption{At t=5}
  \end{subfigure}
  
  \medskip
  \begin{subfigure}[b]{0.33\textwidth}
    \centering
    \includegraphics[width=\linewidth]{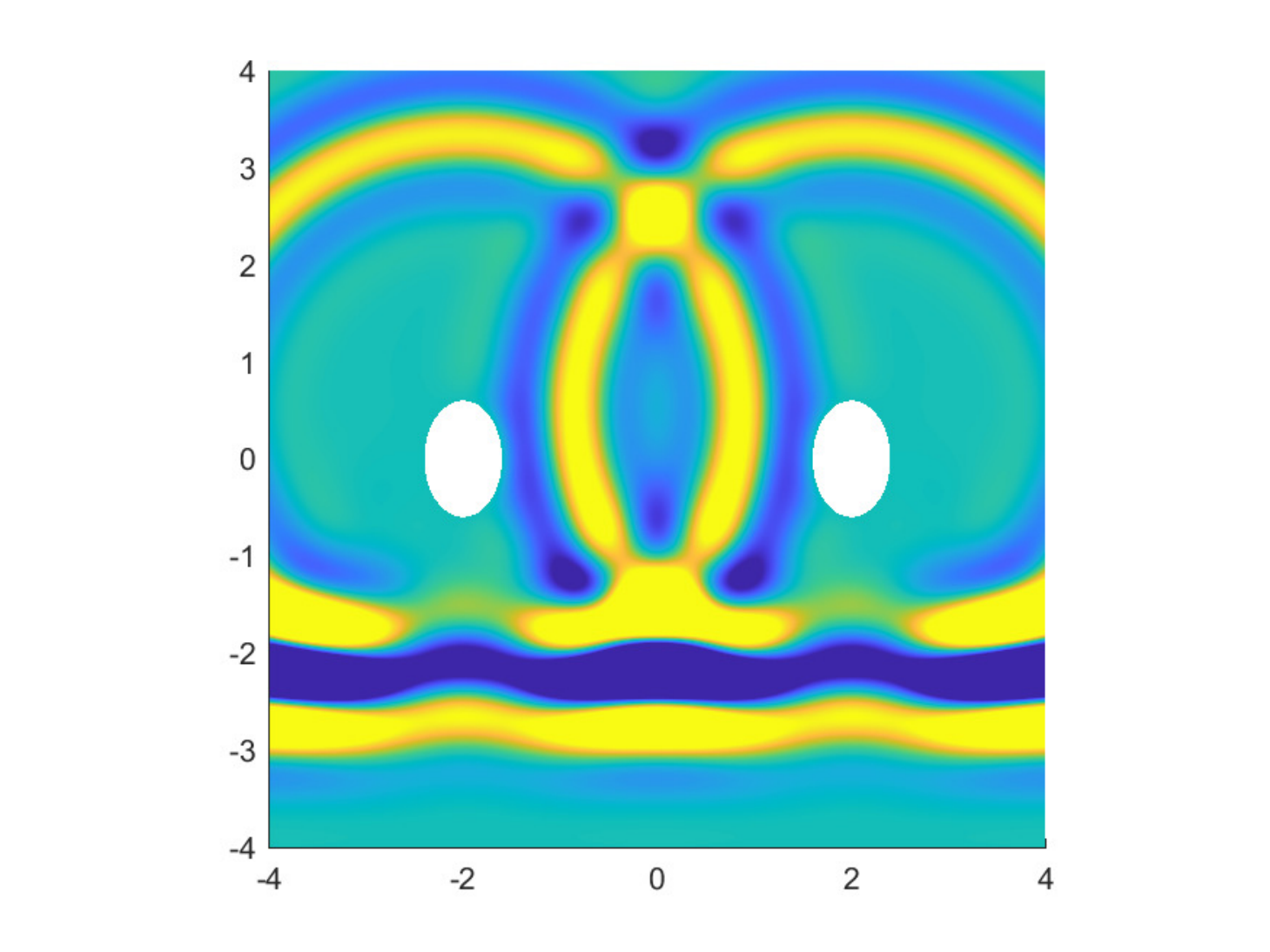}
    \caption{At t=6.25}
  \end{subfigure}\hfil%
  \begin{subfigure}[b]{0.33\textwidth}
    \centering
    \includegraphics[width=\linewidth]{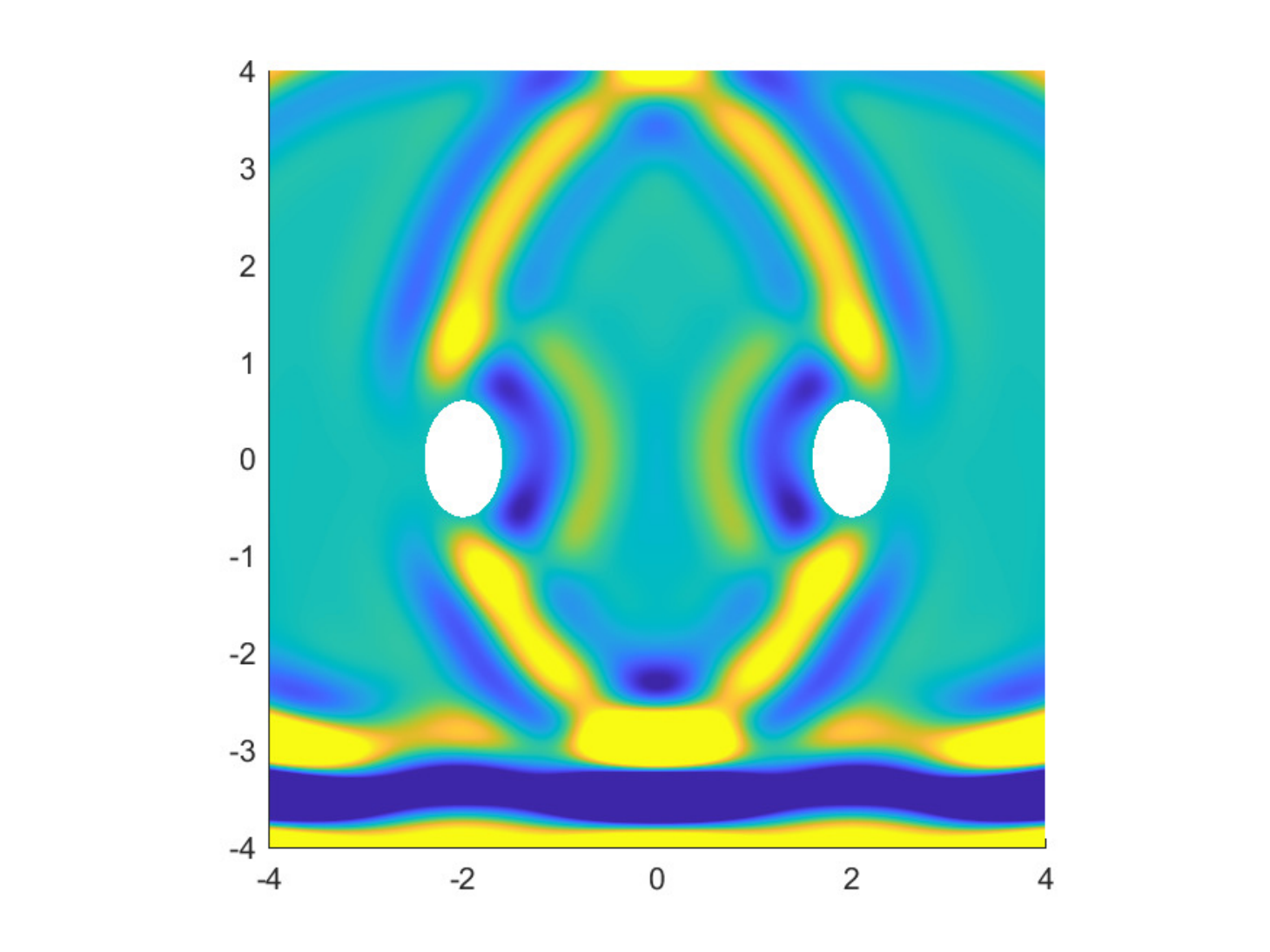}
    \caption{At t=7.5}
  \end{subfigure}\hfil%
  \begin{subfigure}[b]{0.33\textwidth}
    \centering
    \includegraphics[width=\linewidth]{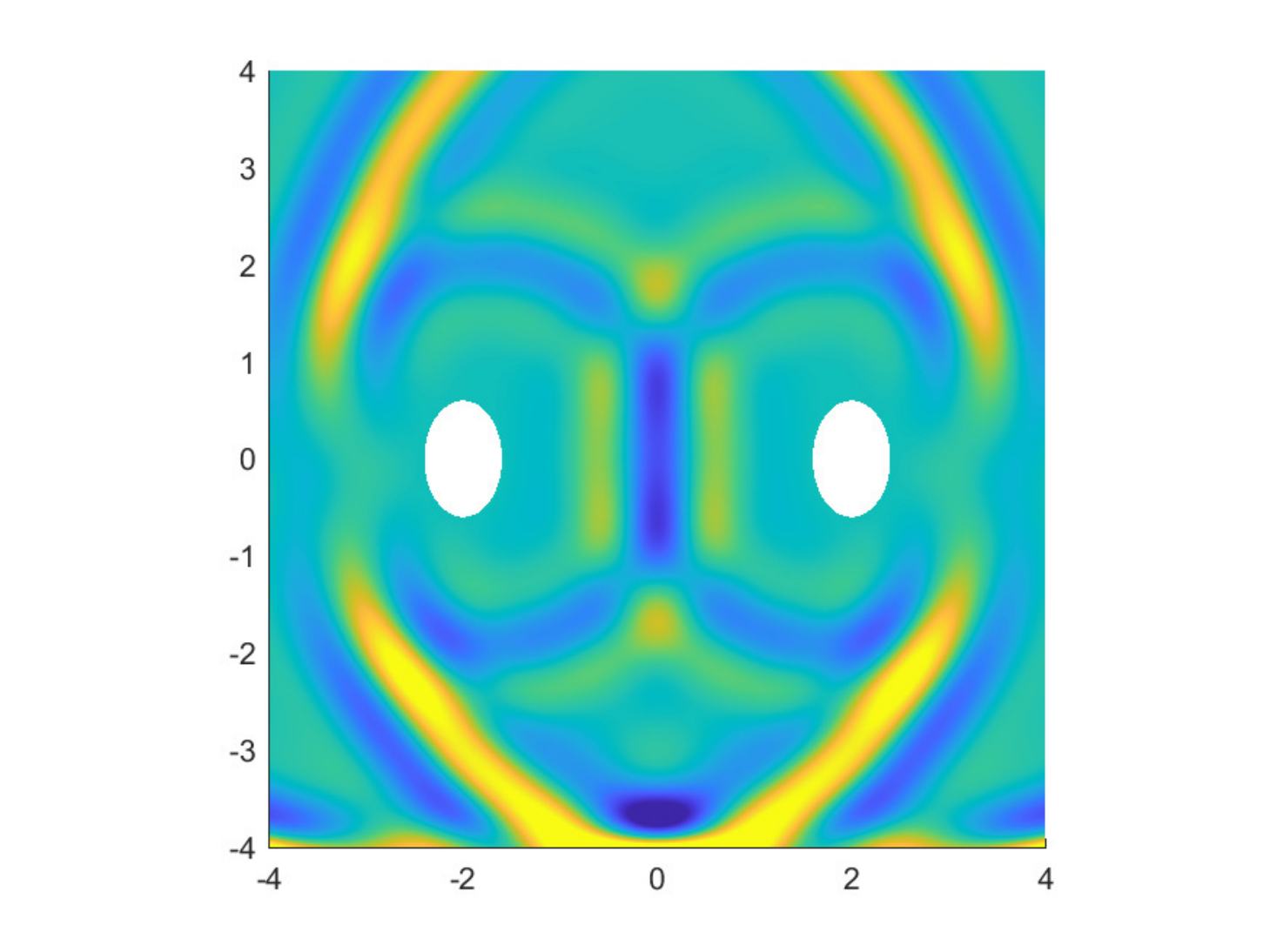}
    \caption{At t=8.75}
  \end{subfigure}%
  \caption{Snapshots of the total field for the scattering by two ellipses.}\label{fig:plot_ellipses}
\end{figure}

Again, we see that a much larger time-step is sufficient to obtain good accuracy when using the modified scheme compared with the standard CQ. However, here we also see a curious behaviour for very large time-steps, where the standard \blue{method, while inaccurate,} also remains reasonably \blue{bounded, whereas} the modified scheme produces very large errors. We have seen this effect often in numerical experiments whenever the time-step is too large. This is not a real deficiency of the method as with such a large time-step the results would anyway be inaccurate but is a property of the modified scheme that one should be aware of.
In Figure~\ref{fig:plot_ellipses} we show snapshots of the total field for the scattering by two ellipses \red{with the frequency $\omega=5$}.

\subsection{Transient wave scattering and the Galerkin method in space}
\begin{figure}
    \centering
    \includegraphics[width=.7\textwidth]{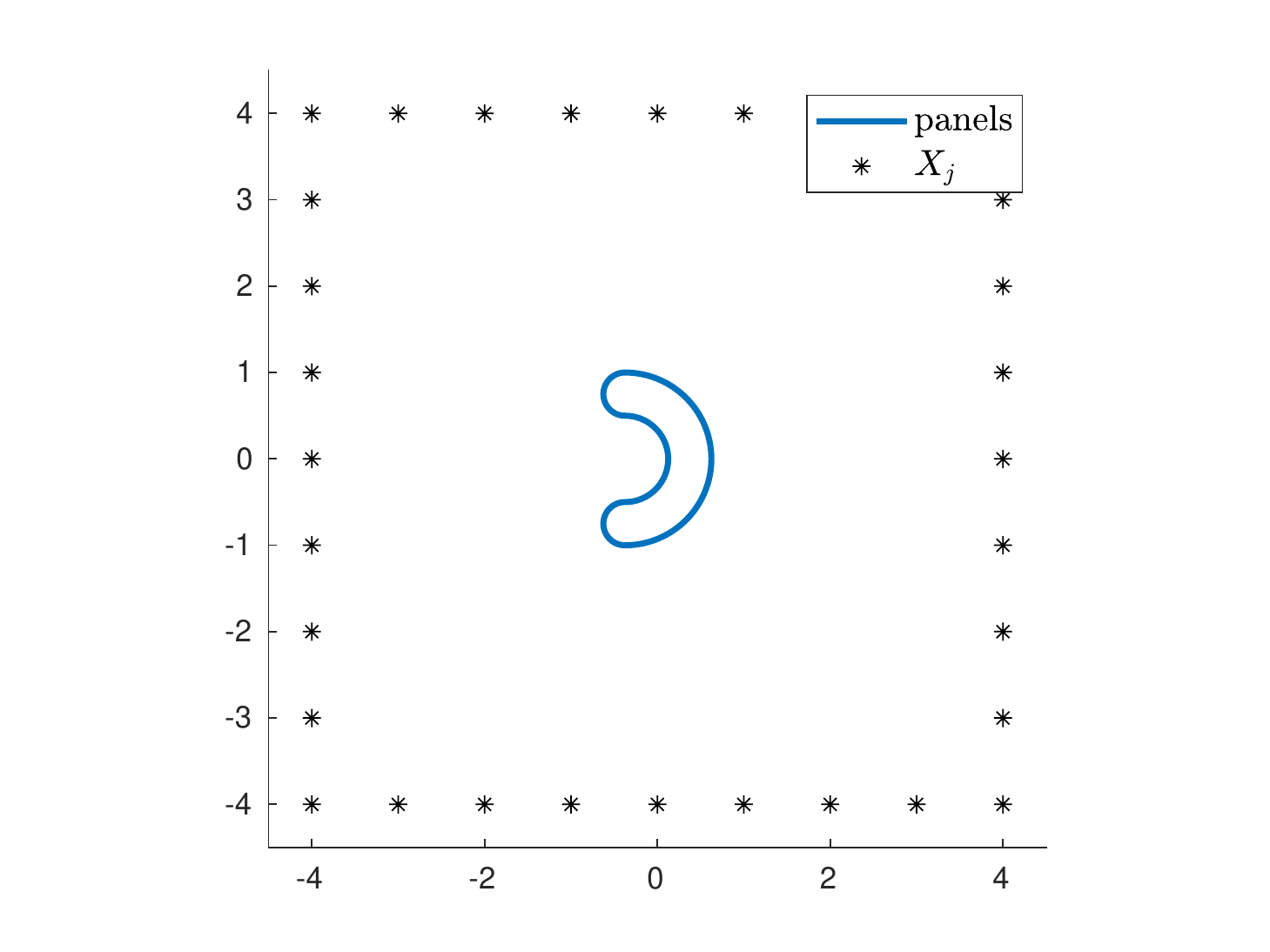}
    \caption{We show the non-convex domain whose boundary is described in \eqref{eq:semicircles} and the points $X_j$ where the error will be computed.}
    \label{fig:semicircles}
\end{figure}

In this section we perform experiments with the Galerkin BEM discretization in space as described in Section~\ref{sec:Gal}, i.e., we use the space of piecewise constant boundary elements. The \blue{time discretization} is as in the previous section and the implementation is based on the formula \eqref{eq:decoupled}. We let the domain $\Omega$ be the non-convex domain seen in Figure~\ref{fig:semicircles}. Its boundary is composed of four semi-circles defined in the complex plane by
\begin{equation}\label{eq:semicircles}
\begin{aligned}
\Gamma_1 &:= \left\{e^{\mi \theta} \;:\; -\frac{\pi}2 \leq\theta< \frac{\pi}2  \right\}\\
\Gamma_2 &:= \left\{\frac14 e^{\mi \theta} + \frac34\mi \;:\; \frac{\pi}2\leq\theta<\frac{3\pi}2 \right\}\\
\Gamma_3 &:= \left\{\frac12 e^{-\mi \theta} \;:\; -\frac{\pi}2 \leq \theta < \frac{\pi}2 \right\}\\
\Gamma_4 &:= \left\{\frac14 e^{\mi \theta} - \frac34\mi \;:\; \frac{\pi}{2}\leq \theta <\frac{3\pi}{2} \right\}.
\end{aligned}
\end{equation}

\begin{figure}%\hspace{-25mm}
\centering
  \begin{subfigure}[b]{0.33\textwidth}
    \centering
    \includegraphics[width=\linewidth]{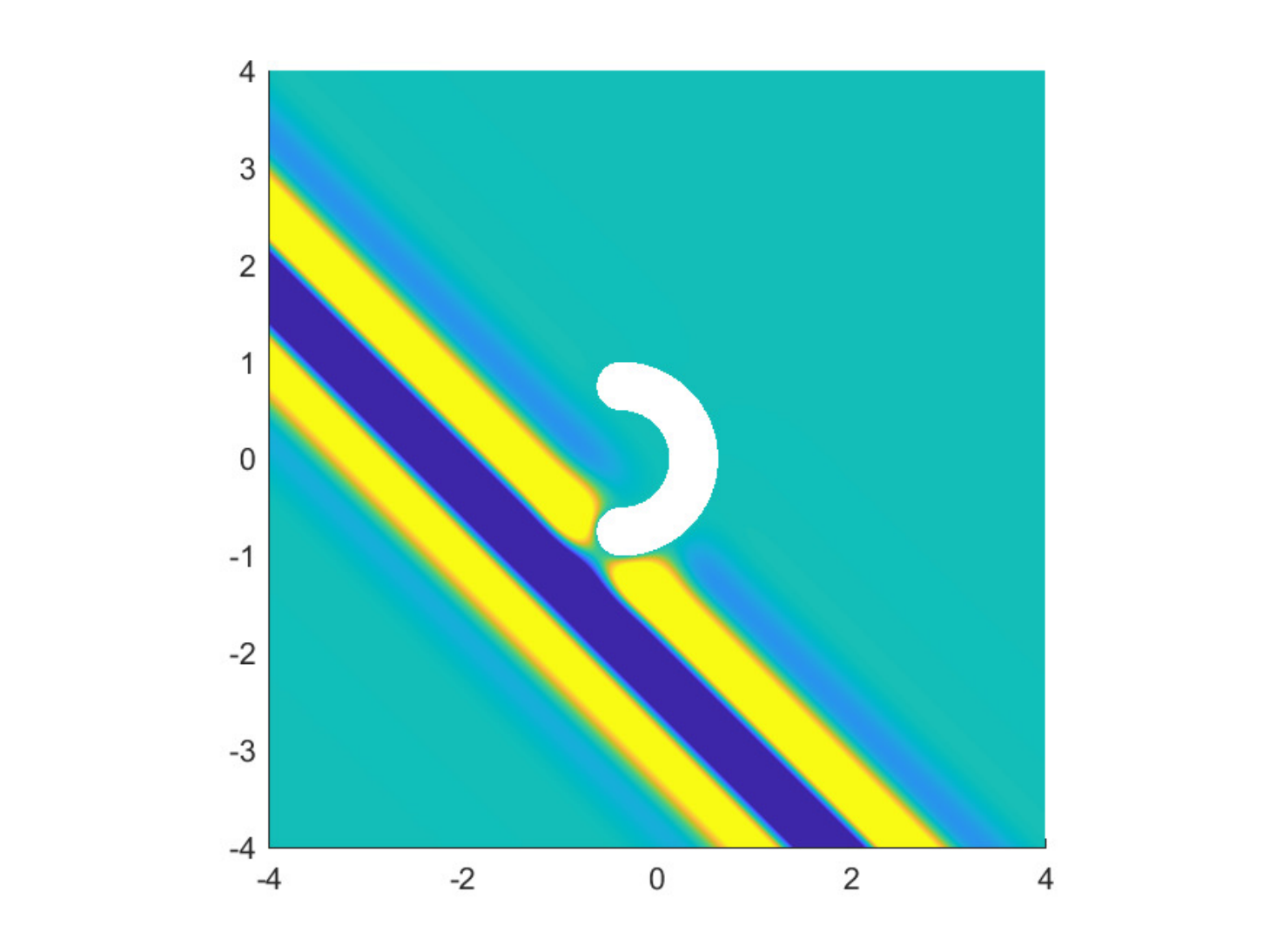}
    \caption{At t=2.5}
  \end{subfigure}\hfil%   
  \begin{subfigure}[b]{0.33\textwidth}
    \centering
    \includegraphics[width=\linewidth]{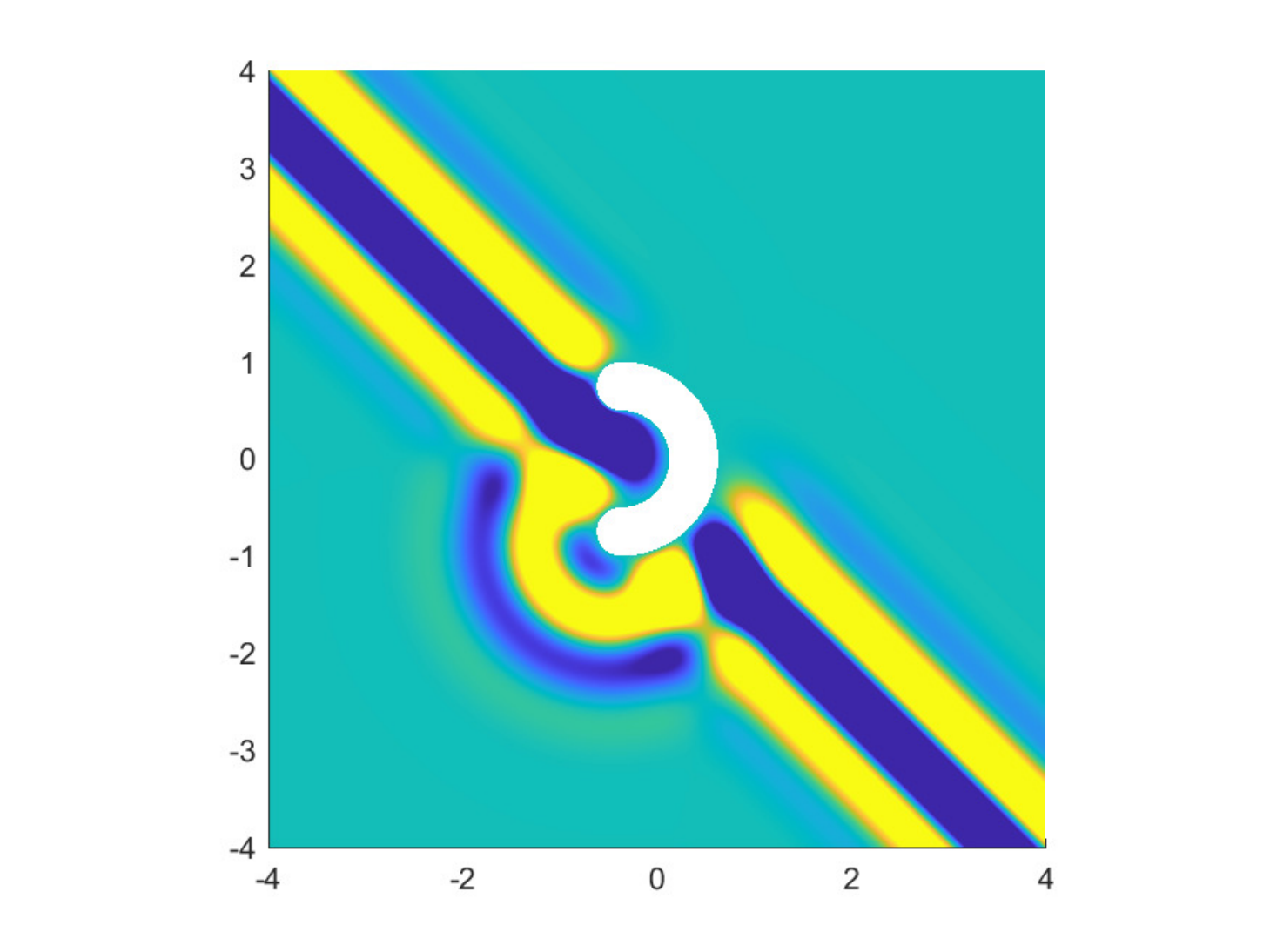}
    \caption{At t=3.75}
  \end{subfigure}\hfil%  
  \begin{subfigure}[b]{0.33\textwidth}
    \centering
    \includegraphics[width=\linewidth]{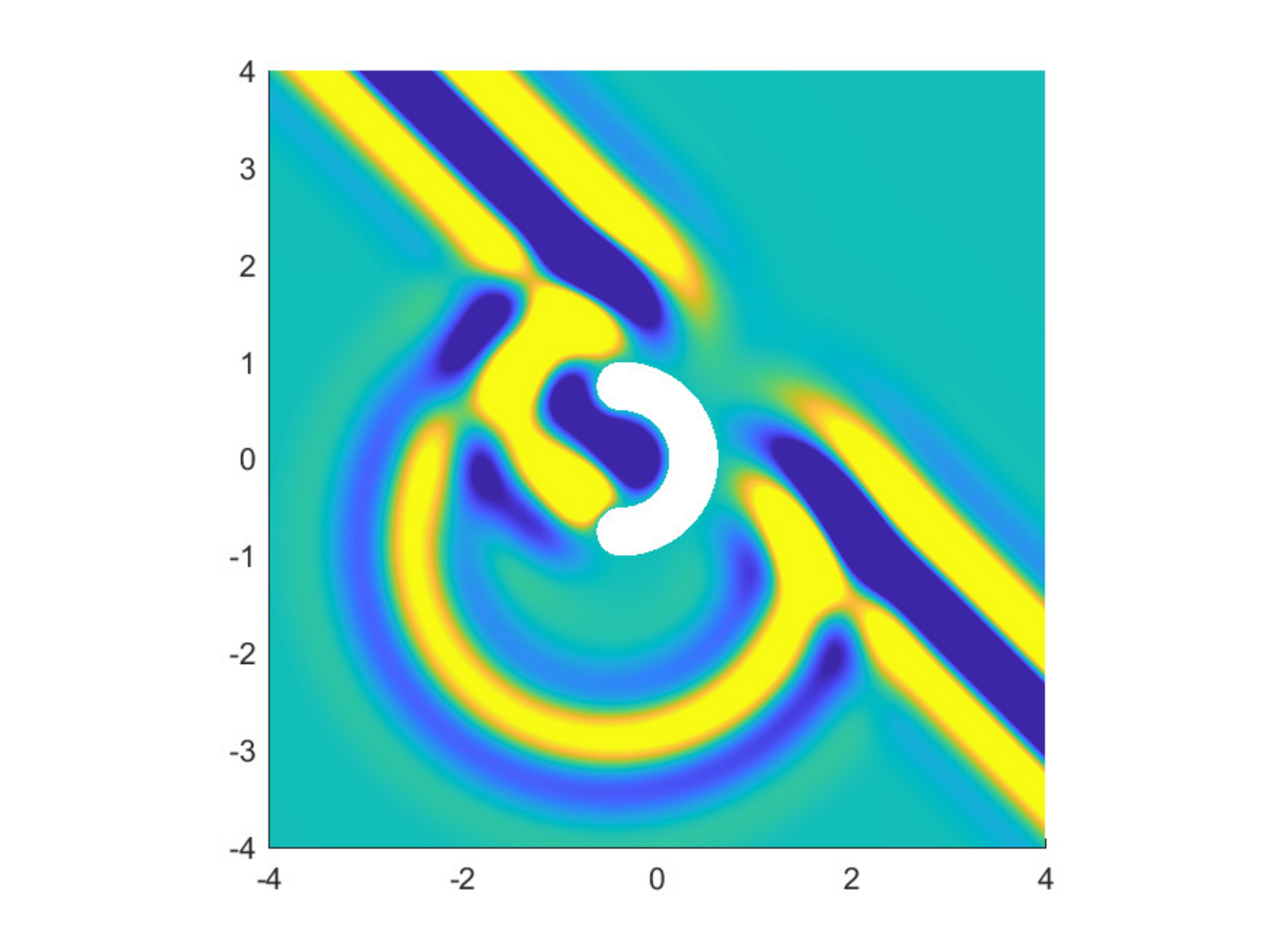}
    \caption{At t=5}
  \end{subfigure}
  
  \medskip
  \begin{subfigure}[b]{0.33\textwidth}
    \centering
    \includegraphics[width=\linewidth]{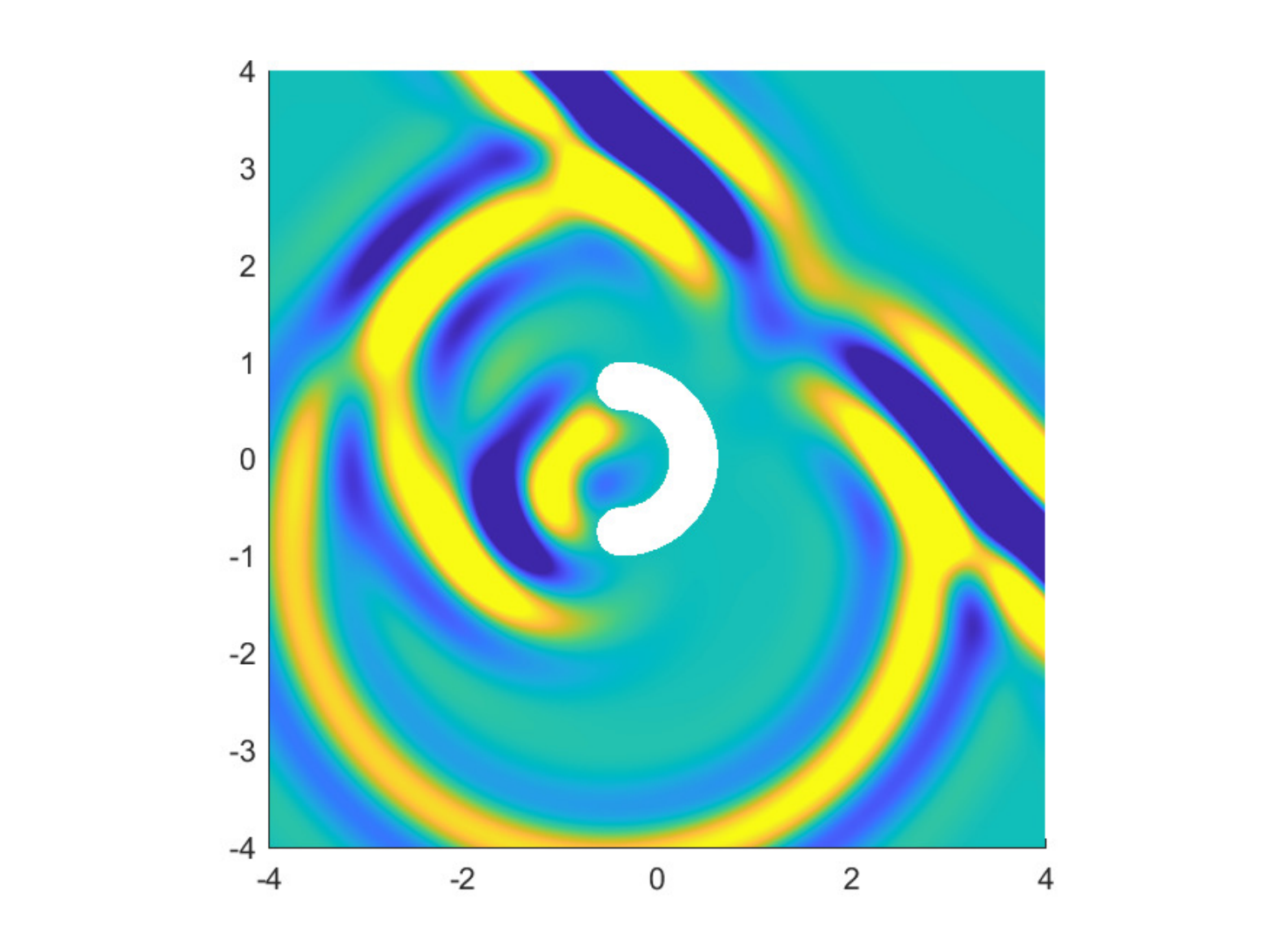}
    \caption{At t=6.25}
  \end{subfigure}\hfil%
  \begin{subfigure}[b]{0.33\textwidth}
    \centering
    \includegraphics[width=\linewidth]{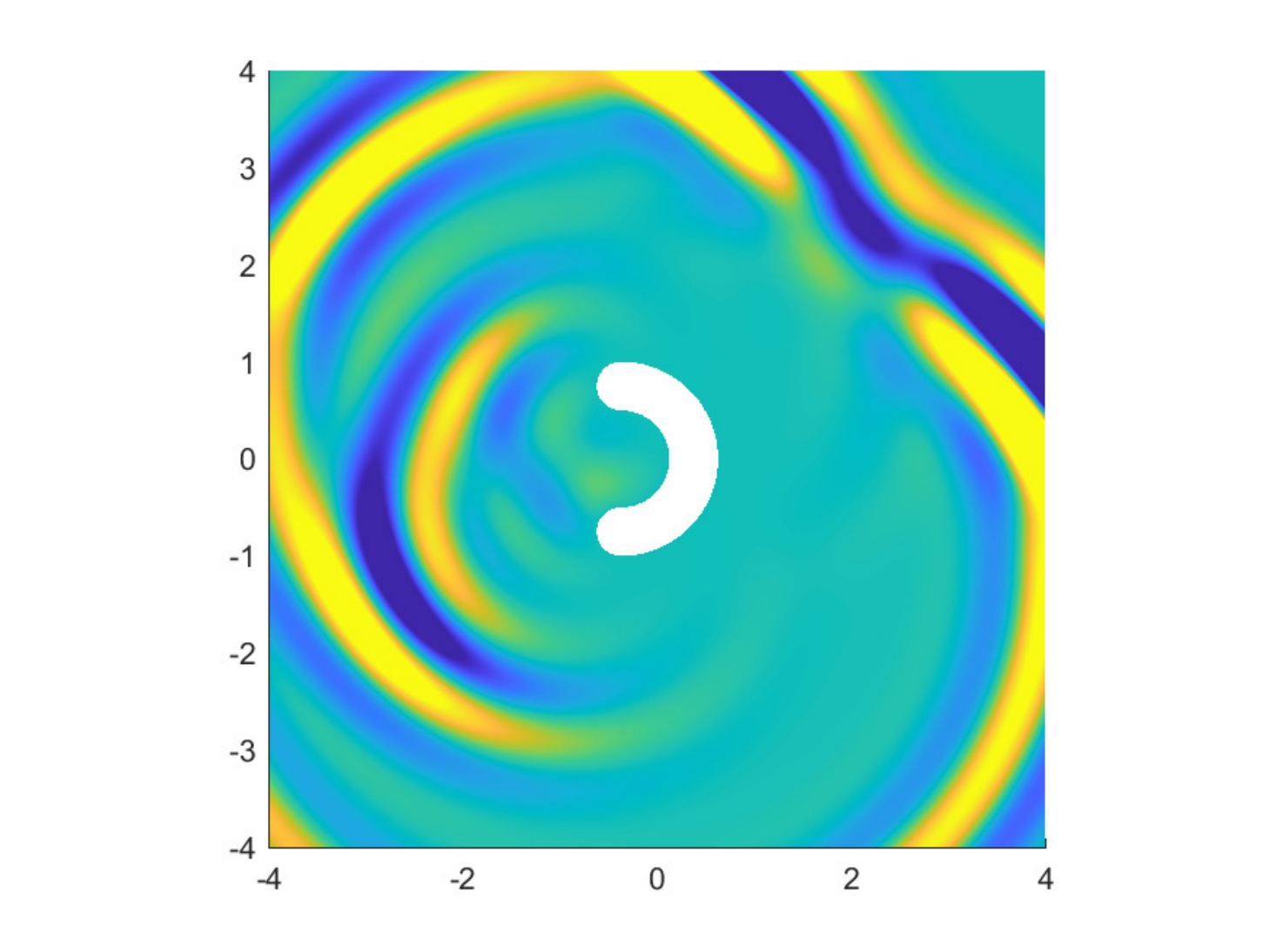}
    \caption{At t=7.5}
  \end{subfigure}\hfil%
  \begin{subfigure}[b]{0.33\textwidth}
    \centering
    \includegraphics[width=\linewidth]{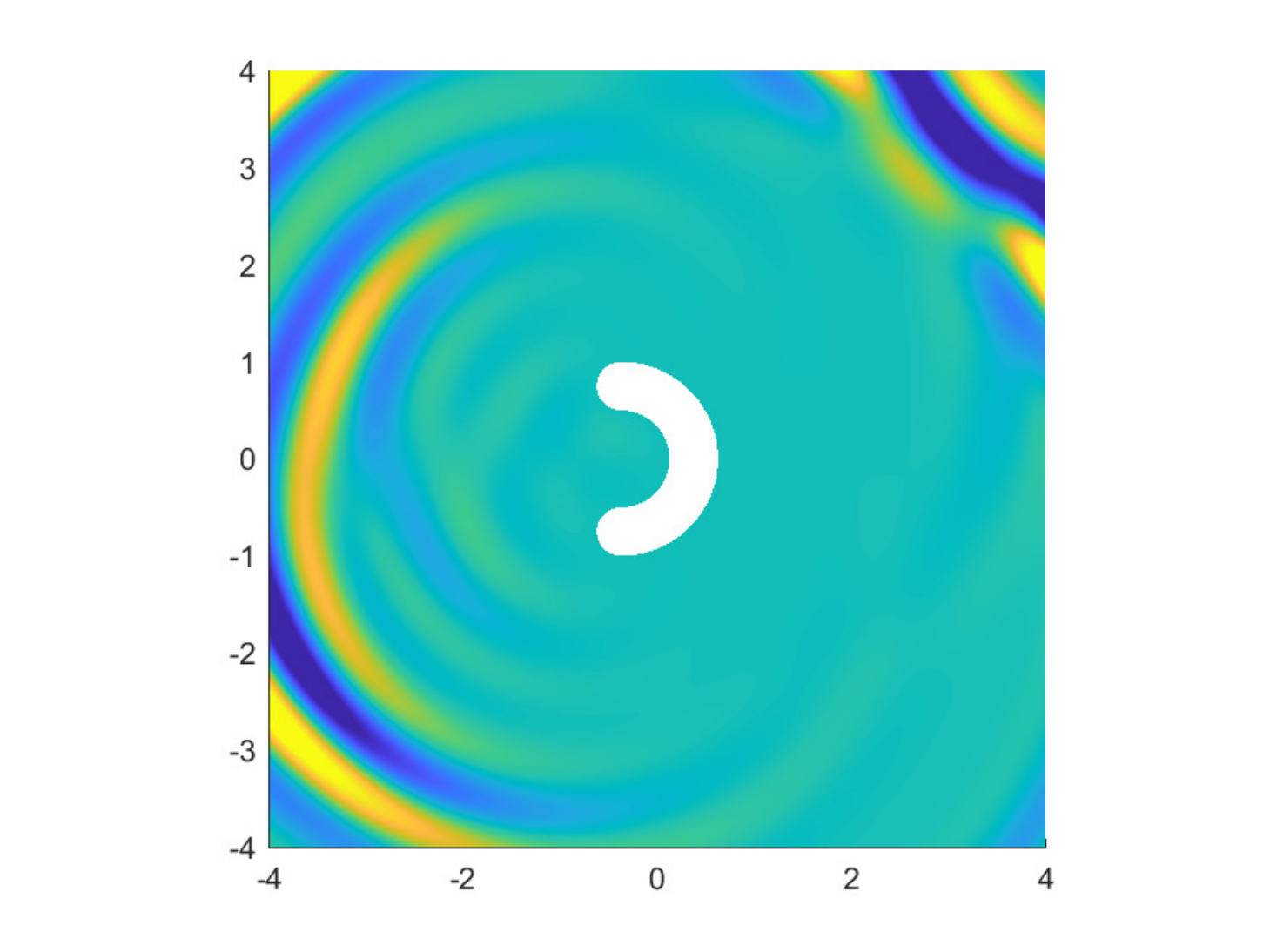}
    \caption{At t=8.75}
  \end{subfigure}%
  \caption{Snapshots of the total field for the scattering by the non-convex domain.}\label{fig:plot_semicircles}
\end{figure}

As the incident wave we use 
\begin{equation}\label{eq:planewave}
\uinc(\red{t,x}) = -\sin(\omega(t-x \cdot \alpha)) g(t-4-x \cdot \alpha)
\end{equation}
with
\[
g(t) = e^{-(t/0.7)^2}, \qquad \alpha = \frac1{\sqrt{2}}(1, 1)^T
\]
and $\omega$ a parameter that we will choose as either $\omega= 1$ or $\omega = 5$. The Dirichlet data is then given by $-\uinc$ and the final time is set to $T = 10$.
Snapshots  of the total field \red{with $\omega=5$} are shown in Figure~\ref{fig:plot_semicircles}.

\begin{figure}
    \centering
    \includegraphics[width=.8\textwidth]{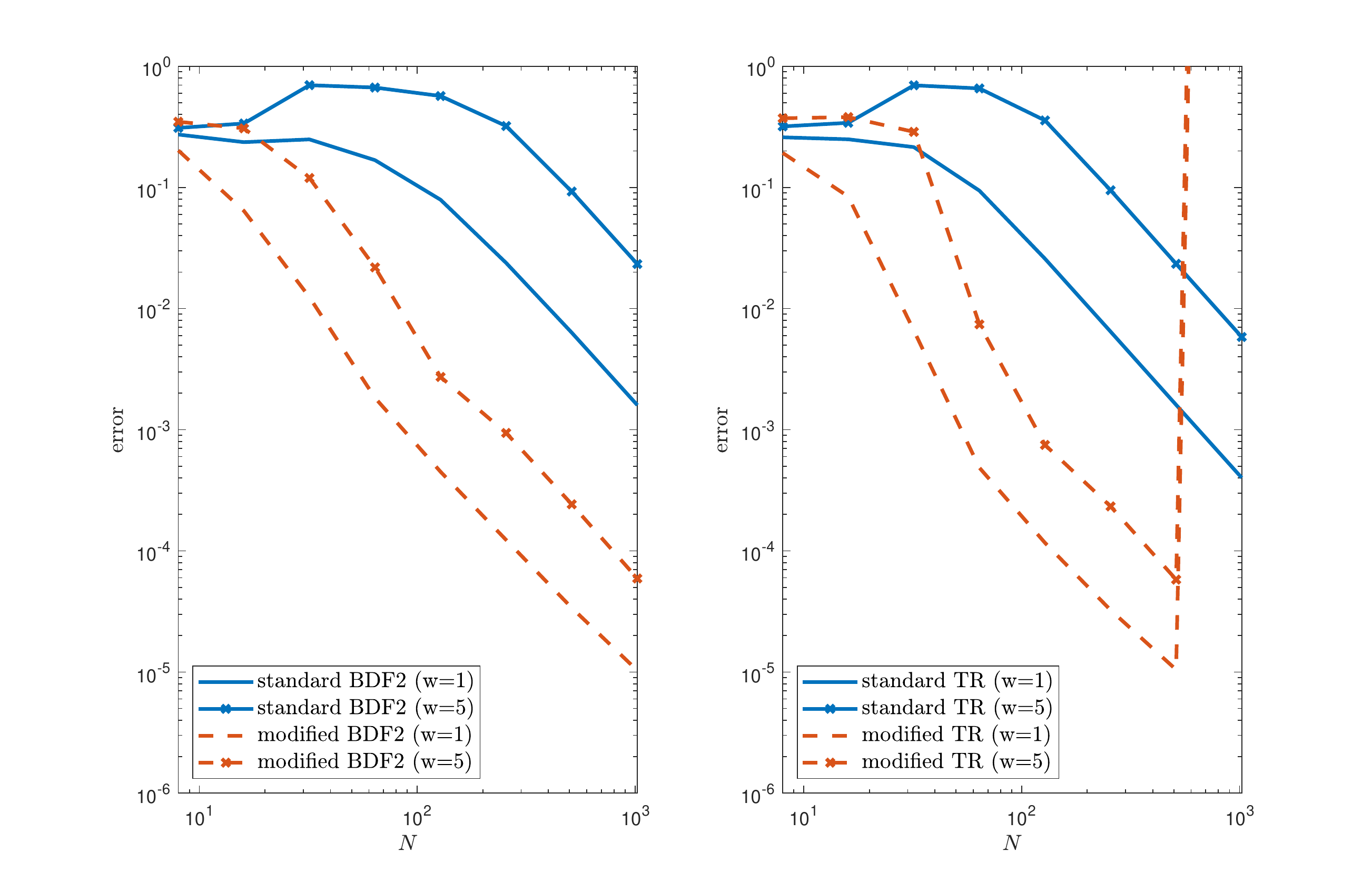}
    \caption{Convergence of the error for Galerkin BEM discretization in space and CQ and modified CQ discretization in time. The scatterer is the non-convex domain in Figure~\ref{fig:semicircles} and the incident wave is the plane-wave \eqref{eq:planewave}. An instability is seen in the modified scheme based on the trapezoidal scheme. This as instability is removed when finer discretization is used in space; see Figure~\ref{fig:conv_semicircles_fine}.}
    \label{fig:conv_semicircles}
\end{figure}

\begin{figure}
    \centering
    \includegraphics[width=.8\textwidth]{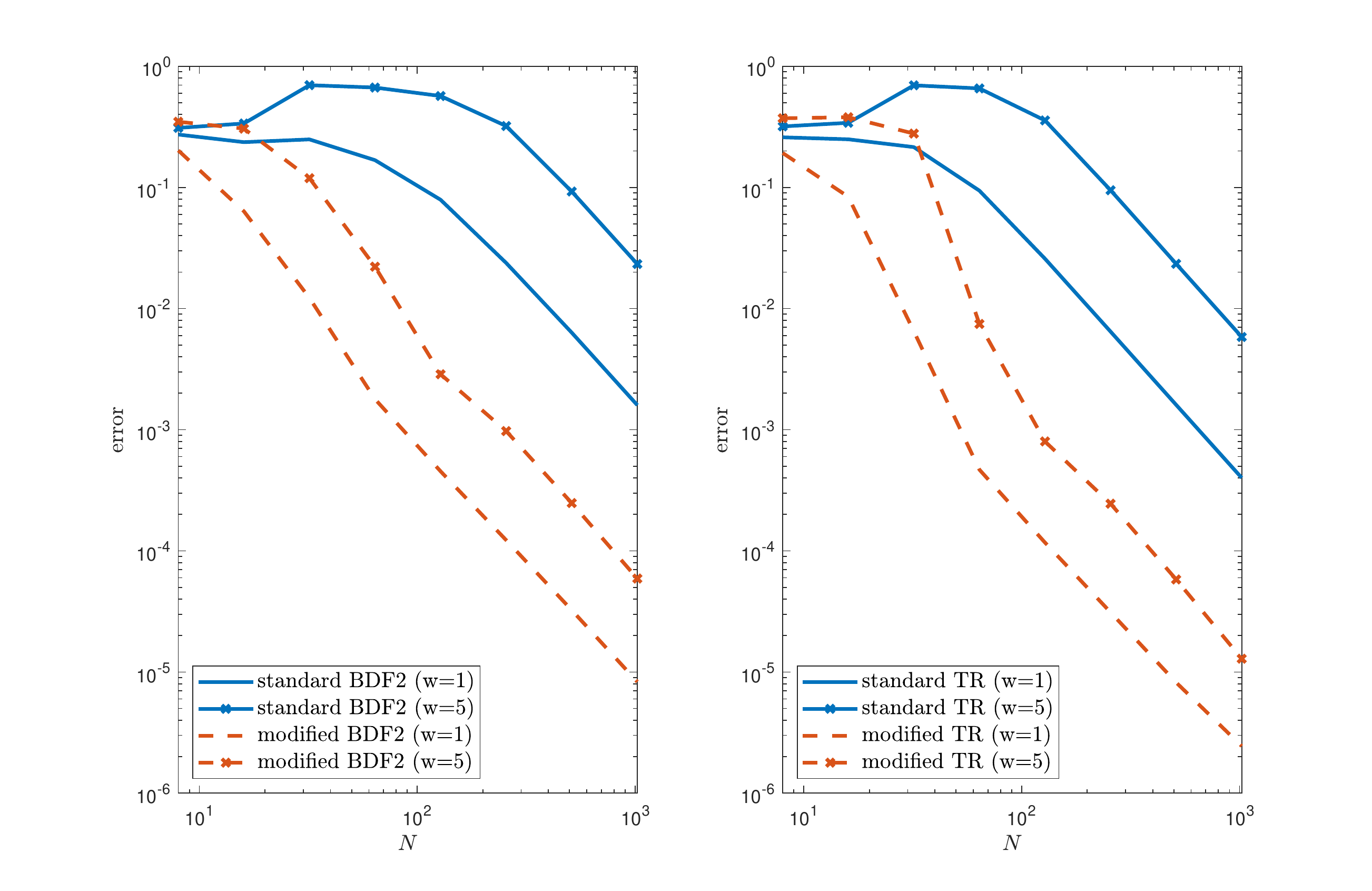}
    \caption{Convergence of the error for Galerkin BEM discretization in space and CQ and modified CQ discretization in time. The scatterer is the non-convex domain in Figure~\ref{fig:semicircles} and the incident wave is the plane-wave \eqref{eq:planewave} with $M = 2000$ and $M_e = 3000$. The finer discretization in space removes the instability seen in Figure~\ref{fig:conv_semicircles}.}
    \label{fig:conv_semicircles_fine}
\end{figure}

The exact solution is not known hence we use a more accurate approximation instead and denote this by $u_{\text{ex}}$. As the error measure we use the maximum error in space at points $X_\ell$, see Figure~\ref{eq:semicircles}, and all time-steps; see \eqref{eq:error}.
For the computation of $u^{\text{ex}}$ a uniform boundary element mesh is used with $M^{\text{ex}} = 1500$ degrees of freedom and $N^{\text{ex}} = 2^{12}$. For the computation of the approximate solution $u^h$ we use again a uniform boundary element mesh with  $M = 1000$ degrees of freedom, whereas the number of time-steps is increased from $N = 2^3$ to $N = 2^{10}$. Again we perform experiments with both BDF2 and trapezoidal based schemes. The convergence results are shown in Figure~\ref{fig:conv_semicircles} and again show considerably better performance of the modified scheme except that for the smallest time-step, the modified scheme based on the trapezoidal scheme becomes unstable. \blue{This is not an unknown problem even for standard convolution quadrature. Its stability is only assured if the quadrature is sufficiently accurate.  } This problem is easily remedied \blue{either by improving the quadrature rule or } by using a finer mesh in space; see Figure~\ref{fig:conv_semicircles_fine} where we increased the spatial degrees of freedom to $M = 2000$ and $M^{\text{ex}} = 3000$. It is by now well-understood that the CQ based on the trapezoidal is more difficult to discretize in space. The reason behind this is that the frequencies
\[
s_\ell = \frac{\delta(\lambda \zeta_{N+1}^{-\ell})}{\tstep}
\]
are of size $\O(\tstep^{-1})$ for the BDF2 scheme and of size $\O(\tstep^{-2})$ for the trapezoidal scheme due to the singularity of $\delta(\zeta)$ at $\zeta = -1$ in the case of the latter; for details see \cite{FJStrap,lbbook}. As the numerical results for the modified BDF2 scheme are comparable to those of the modified trapezoidal scheme, in practice one may prefer to restrict using the BDF2 scheme.

%============= semicircles figures ===============

\subsection{Interior scattering (MFS)}
\begin{figure}
    \centering
    \includegraphics[width=.7\textwidth]{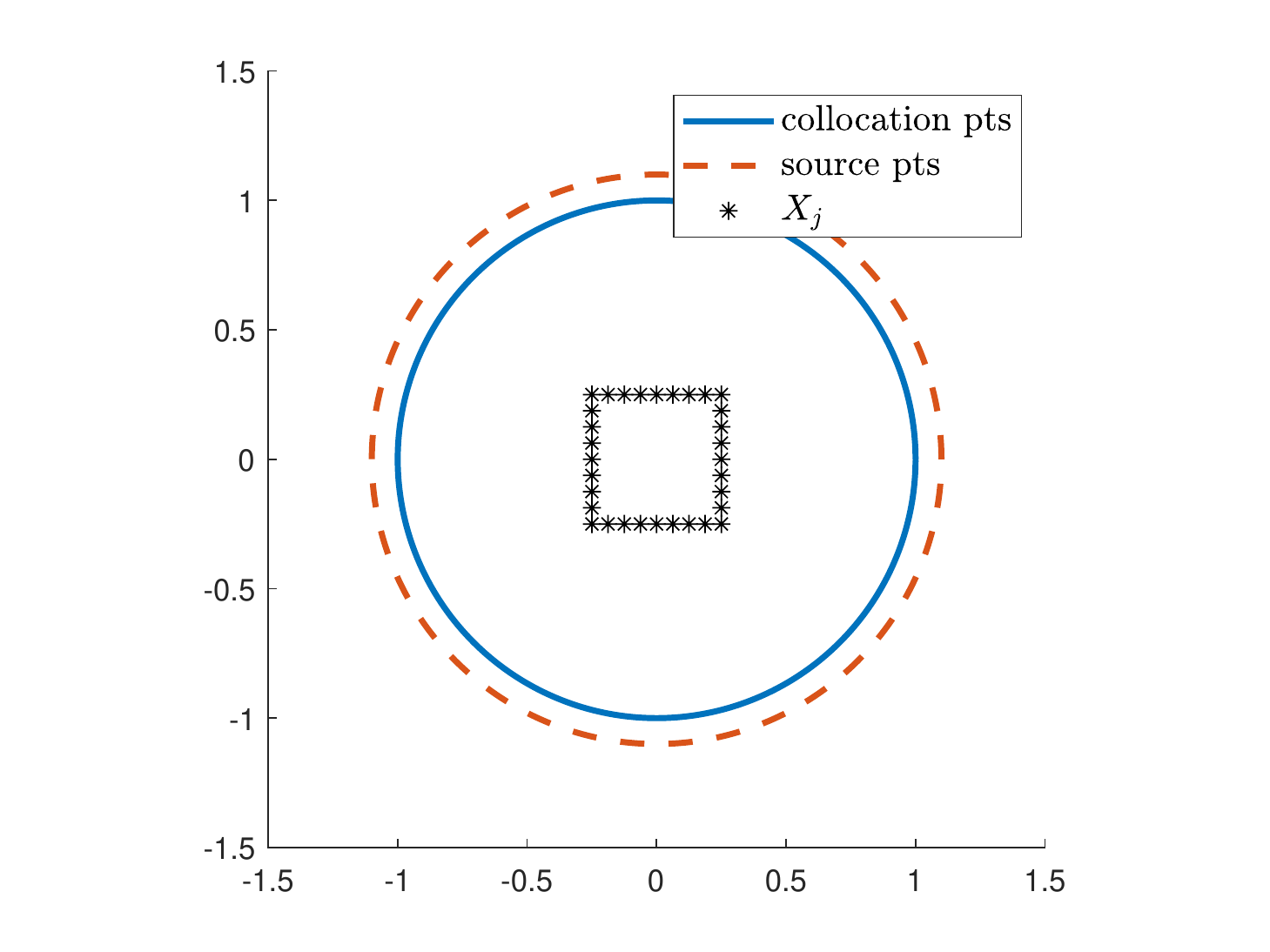}
    \caption{We show the setting for the interior scattering by the unit disk, where the position of the collocation points (solid lines), position of the source points (dashed lines) and the points $X_j$ where we evaluate the error.}
    \label{fig:disk_in_MFS}
\end{figure}

\begin{figure}
    \centering
    \includegraphics[width=.8\textwidth]{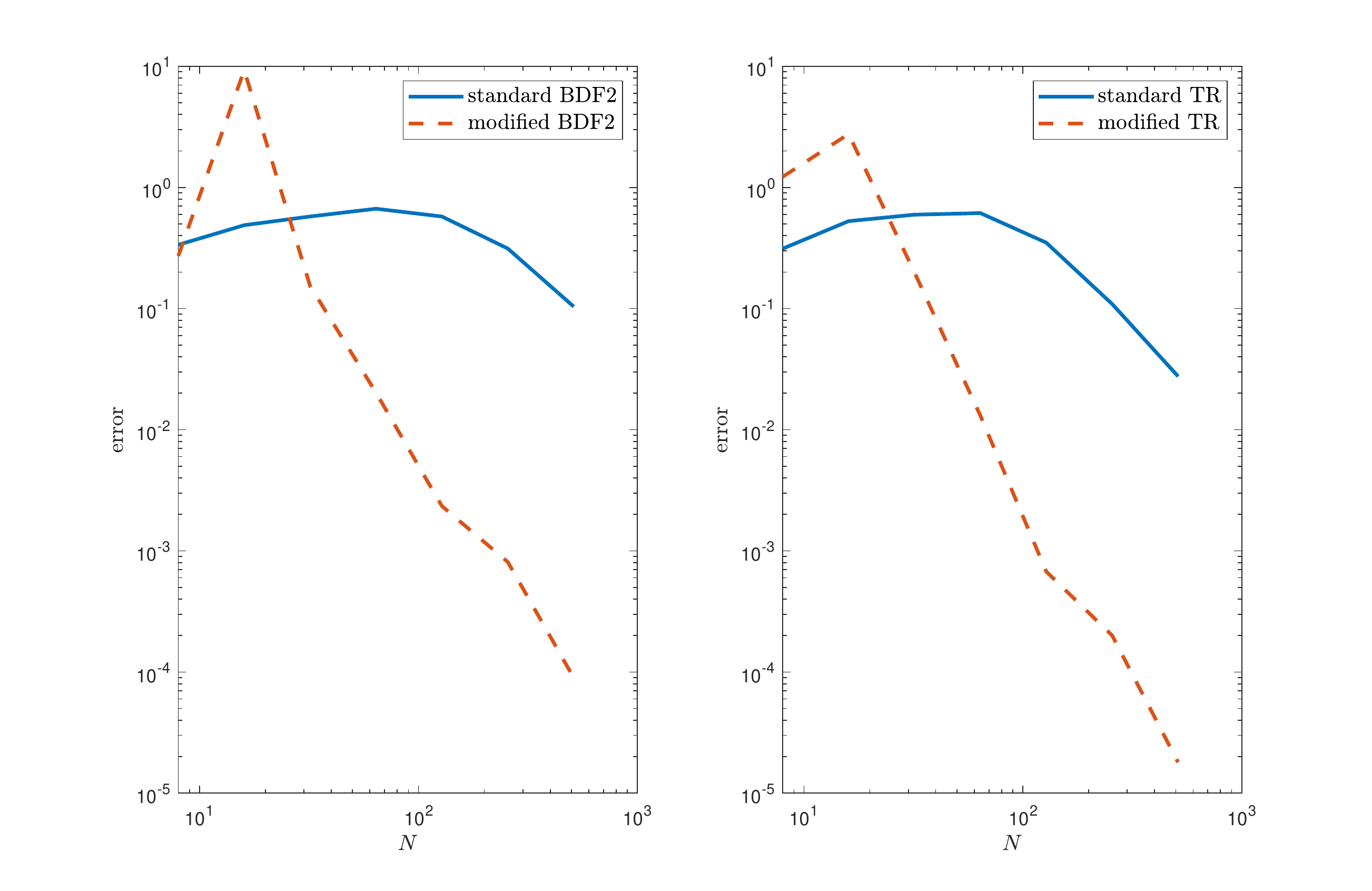}
    \caption{Convergence of the MFS with BDF2 based methods on the left and MFS with trapezoidal rule based methods on the right for the case of interior scattering by the unit disk. }
    \label{fig:conv_disk_in_MFS}
\end{figure}

\blue{
To conclude the numerical experiments, we consider an interior scattering problem. The domain $\Omega$ is the unit disk and the incident wave is given by its initial data
\[
\uinc(0) = e^{-\frac12 a^2 |x-x_0|^2} \qquad \partial_t\uinc(0) = 0,
\]
with
\[
a = 10 \qquad x_0 = (\tfrac14,0).
\]
Note that the choice of parameters is such that $|\uinc(x,0)| \lesssim 10^{-13}$ for $x \in \Gamma = \partial\Omega$, i.e., to almost machine precision the data is initially zero at the boundary. 
Due to the radial form of the initial data, the incident wave can be computed for any $x$, and $t$ using the Hankel transform. Namely, 
\[
\uinc(x,t) = \int_0^\infty F(k) J_0(kr) k \cos(kt) dt \qquad r = |x-x_0|,
\]
where $J_0$ is a Bessel function of the 1st kind and 
\[
F(k) = \frac1{a^2} e^{-\frac{k^2}{2a^2}}
\]
is the Hankel function of the function $e^{-\frac12 a^2 r^2}$.
Thus, the Dirichlet data is given by
\[
g(x,t) = -\uinc(x,t) \qquad x \in \Gamma = \partial \Omega, \; t \geq 0.
\]

We discretize in space using MFS, where 
 use $M = 2000$ collocation points and $K = 1000$ source points. The setting is shown in Figure~\ref{fig:disk_in_MFS}.  To compute the accurate solution $u^{\text{ex}}$ we used $M^{\text{ex}} = 3000$ collocation points and $K^{\text{ex}} = 1500$ source points and in time the modified CQ based on BDF2 with $N^{\text{ex}} = 2^{11}$ time-steps. The convergence of the modified and standard methods is shown in Figure~\ref{fig:conv_disk_in_MFS}.
Again, we see a significant improvement when using the modified method.
Similar results were obtained using the Galerkin method. To save on space we do not show the results here.}

\def\cprime{$'$}

\end{document}